\documentclass[12pt]{article}
\usepackage{amssymb}
\usepackage{amsmath}
\usepackage{amsbsy}
\usepackage{amsthm}
\usepackage{amscd}

\newcommand{\pt}{\partial}

\newcommand{\Pf}{\noindent {\em Proof.} }

\newtheorem{theo}{Theorem}[section]
\newtheorem{rem}[theo]{Remark}

\newtheorem{prop}[theo]{Proposition}
\newtheorem{lem}[theo]{Lemma}
\newtheorem{notat}[theo]{Notation}

\newtheorem{cor}[theo]{Corollary}
\newtheorem{defi}[theo]{Definition}
\newtheorem{defi-lem}[theo]{Definition-Lemma}

\def \C{{\bf C}}
\def \R{{\bf R}}
\def \Z{{\bf Z}}

\def \cZ{\cal Z}
\def \D{{\cal D}}

\def \E{{\cal E}}
\def \F{{\cal F}}

\def \J{{\cal J}}
\def \L{{\cal L}}
\def \K{{\cal K}}
\def \M{{\cal M}}

\def \R{{\bf R}}
\def \cR{{\cal R}}

\def \T{{\cal T}}
\def \V{{\cal V}}

\def \Cst{{$C^*$}}

\def \Lie{{\rm Lie}}
\def \ind{{\rm ind}}

\def \mod{{\rm mod\,}}

\def \Exp{{\rm Exp}}
\def \ext{{\rm ext\,}}
\def \inter{{\rm \,int\,}}
\def\supp{{\rm supp\,}}

\def\hotimes{\hat\otimes}
\def\Ct{Cl_\tau}

\def\c*G{C^*_\Gamma}

\def\l2lf{{L^2_{lf}}}

\def\Kl{{\cal K}_{lc}}
\def\tsig{{\tilde\sigma}}
\def\Cliff{{\rm Cliff}}
\def\dom{{\rm dom}}
\def\const{{\rm const}}

\def\gtc{\mathfrak{c}}
\def\gtd{\mathfrak{d}}
\def\gtK{\mathfrak{K}}
\def\gtk{\mathfrak{k}}

\def\gtN{\mathfrak{N}}
\def\gtn{\mathfrak{n}}

\def\gtg{\mathfrak{g}}
\def\gtG{\mathfrak{G}}
\def\gtD{\mathfrak{D}}
\def\gtd{\mathfrak{d}}

\def\gtR{\mathfrak{R}}

\def\gtr{\mathfrak{r}}

\def\gtG{\mathfrak{G}}

\def\gtZ{\mathfrak{Z}}

\def\gr{\gtg\gtr}
\def\Gr{\gtG\gtr}
\def\GK{\gtG/\gtK}

\def\cqfd{{\vbox{\hrule height 5pt width 5pt}\quad}\penalty -10}

\begin{document}

\title{Coarse pseudo-differential calculus\\ 
and index theory on manifolds\\
with a tangent Lie structure\\
{\small (Version 4;\; June 2025)}}

\date{}

\author{Gennadi Kasparov}

\maketitle

\begin{abstract} We introduce a simplified (coarse) version of pseudo-differential calculus for operators of order zero on complete Riemannian manifolds. This calculus works for the usual H\" ormander $(1,0)$ class operators, as well as for pseudo-differential operators on filtered manifolds. In fact, we develop the coarse PDO calculus on a more general class of manifolds which we call manifolds with a tangent Lie structure. We prove an index theorem for `h-elliptic' operators where the index is not just an integer, but an element of the $K$-homology group of the manifold. 
\end{abstract}

\section{Introduction}

    We present a coarse approach to the pseudo-differential operator (PDO) calculus on complete Riemannian manifolds. This is not a replacement of the usual PDO calculus. It goes along with the usual PDO calculus (\emph{when the usual PDO calculus exists}) and simplifies it in a number of ways. The coarse approach works in the setting of H\" ormander's $(1,0)$ class of operators and van Erp - Yuncken's class of operators on filtered manifolds. In both cases our approach complies with these existing PDO theories.

Roughly speaking, index theory deals with what is called `elliptic operators', which are some kind of bounded (`zero order') operators, invertible modulo some kind of `negative order' operators. In the simplest case, one has a compact manifold $X$ with the usual action of $C(X)$ on $H=L^2(X)$. `Zero order' operators are bounded operators on $H$ which commute with the action of $C(X)$ modulo $\K(H)$ (this is called \emph{pseudolocality}), `negative order' operators are just compact operators. If $X$ is locally compact, we replace $C(X)$ with $C_0(X)$ and keep the pseudolocality condition, but `negative order' operators $T$ are defined as those with the property: $T\cdot C_0(X)\subset \K(H)$ and $C_0(X)\cdot T\subset \K(H)$. Index theory seeks invariants of operators of `zero order' modulo `negative order'. This is very similar to the basic $KK$-theory.

In the classical H\" ormander theory a symbol of a PDO of order $0$ is a bounded matrix-valued function on the (co)tangent manifold $TX$. Its Fourier transform over the tangent direction is a cosymbol -- a continuous field of matrices consisting of bounded multipliers $\{\tsig_x\in \M(C^*(T_x)), x\in X\}$. (Here the tangent space $T_x$ is considered as an abelian group.) In the case of filtered manifolds the groups $T_x$ are replaced with the osculating groups $\gtG_x$. However, there is a smooth family of linear isomorphisms $\{\gtg_x=\Lie(\gtG_x)\simeq T_x\}$. 

According to the coarse PDO approach, to construct an operator on $L^2(X)$ out of a cosymbol $\tsig$ we use a family of local exponential maps (parametrized by $x\in X$). This gives a family of operators on $L^2(X)$. Applying `operator integration' described in section 2 we get an operator on $L^2(X)$ corresponding to the given cosymbol.

There are two basic assumptions on the cosymbol which are required for the coarse PDO construction. The first assumption is the pseudolocality of a cosymbol (i.e. commutation with continuous functions modulo compact operators). The statements given  in section 3 provide sufficient conditions which are easy to use. The second assumption is norm-continuity. This condition is more technical (see definition 5.1). The H\" ormander $\rho,\delta$ class, for $\rho=1, \delta=0$ only, and the van Erp - Yuncken class of operators on filtered manifolds do satisfy this condition (see section 6).

\medskip 

We present the coarse PDO approach in the framework of manifolds which we call \emph{manifolds with a tangent Lie structure}. This is a generalization of filtered manifolds.

The simplest example of such manifold $X$ is any solvable Lie group $\gtG$. The tangent space $T_x$ at any point $x\in X$ is isomorphic to the Lie algebra $\gtg$ of the group $\gtG$. The family of these isomorphisms $T_x\simeq \gtg$ is smooth in $x$. A concrete example of a possible full PDO calculus on such manifold which is not a H\" ormander calculus is indicated in subsection 6.4.

In general, a tangent Lie structure on a manifold $X$ is defined as a smooth field of Lie algebras $\{\gtg_x, x\in X\}$ assigned to points of $X$ such that \emph{linearly} $\gtg_x\simeq T_x(X)$ for any $x\in X$. We can denote by $\{\gtG_x\}$ a smooth bundle of simply connected Lie groups with $\Lie(\gtG_x)=\gtg_x$. In the filtered manifold case $\gtG_x$ are called osculating groups. 

For a differential operator $\D$ on $X$ which has first order in the usual sense, by freezing coefficients of $\D$ at each point $x\in X$ we get a smooth family of constant coefficient differential expressions on the field of spaces $\{T_x(X)\}$. Since $\gtg_x\simeq T_x(X)$, we can interpret this family as a smooth family of elements of the field $\{\gtg_x\}$ or as a family of first order left-invariant differential operators $\{D_x\}$ on the field $\{\gtG_x\}$ and call this family a cosymbol of $\D$. 

However, a full differential calculus is unlikely to exist in this generality without additional assumptions on the field $\{\gtg_x\}$. One needs the property that the space of `operators of order $k$ modulo lower order' is isomorphic to the space of `symbols of order $k$ modulo lower order'. This is difficult to imagine when even the definition of the `order' is unclear. Therefore in general we cannot discuss the full PDO calculus for manifolds with a tangent Lie structure yet. \emph{On the other hand, the coarse PDO calculus does exist in full generality on manifolds with a tangent Lie structure.}

To be more precise about the definition of a tangent Lie structure,
let $\gtK_x$ be a maximal compact subgroup in $\gtG_x$. We want a cosymbol of an operator of order $0$ in our calculus to be a multiplier of $C_r^*(\gtG_x)$ at each point $x\in X$. The space $L^2(\gtG_x)$ splits into a direct sum $\oplus_\pi L^2(\gtG_x\times_{\gtK_x} V_\pi)$ according to the decomposition of $L^2(\gtK_x)$ into irreducible representations. The algebra $C^*(\gtG_x)$ splits accordingly. We simplify the picture by replacing $L^2(\gtK_x)$ with a finite sum of $V_\pi$'s. So we define the tangent Lie structure not by setting $\Lie(\gtG_x)\simeq T_x(X)$, but by setting $\Lie(\gtG_x)/\Lie(\gtK_x)\simeq T_x(X)$. 

In a geometric sense, because $\gtG_x/\gtK_x$ is diffeomorphic to a Euclidean space (cf. \cite{Ab}), this choice is reasonable. Of course, for solvable Lie algebras $\gtg_x$, we still have $\gtg_x\simeq T_x(X)$ because $\gtK_x=\{1\}$. On the other hand, for semisimple $\gtg_x$, the symbol theory of differential operators in this approach should be related with $\gtG_x$-invariant differential operators on symmetric spaces $\gtG_x/\gtK_x$ and with spherical functions (see remark 3.9 in section 3 and \cite{Helg}, ch. 10).

As we mentioned above, index theory on manifolds with a tangent Lie structure is developed in the present paper based on the coarse PDO approach. Differential operators and their parametrices will be only briefly discussed in section 9 in relation with the index theorem. But essentially we work only with pseudo-differential operators of order $0$ and of negative order, and we prove an index theorem which generalizes the usual $K$-theoretic Atiyah-Singer index theorem for elliptic operators in the usual form \cite{Ka16}, section 4.

Our index formula states that the index class of an h-elliptic operator of order $0$, as an element of the $K$-homology group of the manifold, is equal to the $KK$-product of its cosymbol class by the Dolbeault element. The Dolbeault element is not the classical one but is a natural generalization of it based on the generalized Connes -Thom isomorphism of \cite{Ka88}, section 5. The combination of the $KK$-theoretic approach and the coarse index theory approach allows a very simple proof of this index formula, saving a lot of effort of the previous work (\cite{vErp1-2,Moh20,Moh22}). 

\medskip
 
There are positive and negative sides in the coarse approach. On the positive side, the whole PDO theory gets much simpler. This approach is naturally related with $KK$-theory and straightforwardly leads to an index theorem. On the negative side, this approach is not convenient in treating differential operators and their parametrices. Therefore, the best option is to have the full and coarse approaches available together, as in the case of H\" ormander's and van Erp - Yuncken's PDO calculi.

The paper is organized as follows. Section 2 contains technical results on operator integration used in the coarse PDO calculus, including the group averaging of operator functions. Section 3 contains technical results on the pseudolocality property. In section 4 we introduce manifolds with a tangent Lie structure and some geometric technical tools. Section 5 contains the coarse PDO construction. In section 6, we discuss the relation of our coarse approach with the H\" ormander $\rho=1, \delta=0$ and the van Erp - Yuncken calculi. Section 7 provides $K$-theoretic preliminaries for the index theorems. Index theorems are proved in section 8. Section 9 extends the main index theorem to differential operators on filtered manifolds.

\section{Operator techniques}

This section contains an operator integration technique for the construction of pseudo-differential operators. We also discuss some background for the definition of cosymbols.

Throughout this section, $X$ will be a second countable, locally compact, $\sigma$-compact space, $D$ a \Cst-algebra, and $\phi:C_0(X)\to \M(D)$ a homomorphism such that $\phi(C_0(X))\cdot D$ is dense in $D$. This implies that $\phi$ extends to a unital homomorphism $\phi:C_b(X)\to \M(D)$, where $C_b(X)$ is the \Cst-algebra of all bounded continuous functions on $X$. We will denote the set of all compactly supported continuous functions on $X$ by $C_c(X)$.

\subsection{Support of an operator}

\begin{defi} The support of an element $F\in \M(D)$ is the smallest closed subset of $X\times X$, denoted $\supp(F)$, such that for any $a,b\in C_c(X)$, one has $\phi(a)F\phi(b)=0$ as soon as $\supp(a)\times \supp(b)\cap \supp(F)=\emptyset$. If $\supp(F)$ is compact, $F$ will be called compactly supported. An element $F\in \M(D)$ will be called properly supported if both projections $p_1:\supp(F)\to X$ and $p_2:\supp(F)\to X$ are proper maps.

Similar definitions also hold (and will be most often used below) for the support of $F\,\mod D\in \M(D)/D$.
 \end{defi}
 
 \begin{rem} {\rm Suppose $F$ is properly supported and $a\in C_c(X)$. Then both $\phi(a)F$ and $F\phi(a)$ are compactly supported. Indeed, it is easy to check that $\supp(\phi(a)F)\subset p_1^{-1}(\supp(a))\cap \supp(F)$, which is compact. Similarly for $F\phi(a)$.
 Actually one can define properly supported $F$ by the condition: for any $a\in C_c(X)$ there is $b\in C_c(X)$ such that $F\phi(a)=\phi(b)F\phi(a)$ and $\phi(a)F=\phi(a)F\phi(b)$.}
\end{rem}

We will need the following fact which was stated without proof in \cite{Ka75}, section 3, proposition 4.

\begin{theo} An element $F\in \M(D)/D$ commutes with $C_c(X)$ (and hence with all $C_b(X)$) if and only if the support of $F$ in $\M(D)/D$ belongs to the diagonal $\Delta$ of $X\times X$.
\end{theo}

\Pf For any $a,b\in C_c(X)$, $\supp(a)\times \supp(b)\cap \Delta=\supp(a)\cap \supp(b)$ (when $\Delta$ is identified with $X$). So the assumption that $\supp(F)\subset \Delta$ means that $\phi(a)F\phi(b)=0$ when $\supp(a)\cap\supp(b)=\emptyset$. Hence to prove the `only if' part we have to show that when $F$ commutes with $C_c(X)$, then for any $a,b\in C_c(X)$ with $\supp(a)\cap\supp(b)=\emptyset$, one has: $\phi(a)F\phi(b)=0$. This is obviously true. 

For the `if' part, we will take a faithful representation of $\M(D)/D$ in a Hilbert space $H$. We denote the composition $C_b(X)\to \M(D)/D\to \L(H)$ by $\psi$. Obviously, $\psi$ extends to characteristic functions of Borel sets in $X$. For simplicity, we keep the notation $F$ for the image of $F$ in $\L(H)$. Also for a Borel set $U\subset X$, we will denote the image by $\psi$ of the characteristic function of $U$ in $\L(H)$ by $P(U)$. The closure of the set $U$ will be denoted by $\bar U$.

The condition that for any $a,b\in C_c(X)$ with $\supp(a)\cap\supp(b)=\emptyset$, one has: $\phi(a)F\phi(b)=0$ implies that if $U,V\subset X$, and the intersection of the closures $\bar U\cap \bar V$ is empty, then $P(U)FP(V)=0$ in $\L(H)$. 

Let $a\in C_c(X)$. Without loss of generality, we can assume that $0\le a\le 1$. To prove that $\psi(a)F=F\psi(a)$, we will first approximate $a$ by a step-function. Fix a positive integer $n$. Let $U_i=\{x\in X|a(x)<(2i+1)/2n\}$. We have: $U_0\subset U_1\subset...\subset U_n=X$. Note that for $|i-k|\ge 2$, $\overline{U_{i+1}-U_i}\cap\overline{U_{k+1}-U_k}=\emptyset$. Since for any $x,y\in \overline{U_{i+1}-U_i}$ one has: $|a(x)-a(y)|\le 1/n$, we get: 
$||\psi(a)-\sum_{i=0}^n P(U_i-U_{i-1})\cdot (i/n)||\le 1/n$ in $\L(H)$.

We rewrite the sum in the last expression as $\sum_{i=0}^n i/n[P(U_i)-P(U_{i-1})]=1/n[\sum_{i=0}^niP(U_i)-\sum_{i=0}^{n-1}(i+1)P(U_i)]=((n+1)/n)-1/n\sum_{i=0}^nP(U_i)$ because $P(U_n)=1$. To prove the commutation property $[\psi(a),F]=0$, it is enough to show that the norm of the commutator $[\sum_{i=0}^n P(U_i),F]$ is bounded in $\L(H)$ by a constant which does not depend on $n$. 

We have: $P(U_i)F-FP(U_i)=P(U_i)FP(X-U_i)-P(X-U_i)FP(U_i)$. It is enough to evaluate $||\sum_{i=0}^n P(U_i)FP(X-U_i)||$ (the other sum evaluates similarly). Note that the sum actually goes to $n-1$ because $X-U_n=\emptyset$. We can replace in this sum $P(U_i)$ with $P(U_i-U_{i-1})$ because $\bar U_{i-1}\cap\overline{X-U_i}=\emptyset$. Also we can replace $P(X-U_i)$ with $P(U_{i+1}-U_i)$ for a similar reason. Then the sum $||\sum_{i=0}^{n-1} P(U_i)FP(X-U_i)||$ will become $||\sum_{i=0}^{n-1} P(U_i-U_{i-1})FP(U_{i+1}-U_i)||$. 

Let us denote $P(U_i-U_{i-1})$ by $P_i$ (we consider $U_{-1}=\emptyset$). We have: $\sum_{i=0}^nP_i=1, P_i^*=P_i, P_i^2=1, P_iP_j=0$ for $i\ne j$. We need to evaluate the norm of $S=\sum_{i=0}^{n-1} P_iFP_{i+1}$. We have: $S^*S=\sum_{i,j=0}^{n-1} P_{i+1}F^*P_iP_jFP_{j+1}\le \sum_{i=1}^n P_iF^*FP_i$. The norm of the last expression is $\le ||F||^2$. \cqfd

\subsection{Riemann operator integration}

This kind of operator integration was essentially the subject of section 3 of \cite{Ka75}. We will add here a few more details.

\begin{defi} We will call an element $F\in \M(D)$ locally compact if for any $a\in C_0(X)$, both $\phi(a)F$ and $F\phi(a)$ belong to $D$. The set of locally compact elements will be denoted $D_{lc}$. In the case of $D=\K(H)$ for a Hilbert space $H$, the notation for locally compact elements will be $\K_{lc}(H)$.

We denote by $Q_{C_0(X)}(D)$ the subalgebra of $\M(D)$ consisting of elements $T\in \M(D)$ which commute with $\phi(C_0(X))$ modulo $D$. The algebra $D_{lc}$ is a two-sided ideal in $Q_{C_0(X)}(D)$.
\end{defi}

In the assumptions described at the beginning of this section, let $F:X\to \M(D)$ be a bounded norm-continuous map such that $F(x)$ commutes with $\phi(C_0(X))$ modulo $D$. We will construct a Riemann type operator integral $\int_X F(x)d\phi\in Q_{C_0(X)}(D)/D_{lc}$ which has the following properties (cf. \cite{Ka75}, section 3, theorem 1):

\begin{theo} $1^{\rm o}$. If $||F(x)|| \le c$ for all $x\in X$ then $\int_X F(x)d\phi \le c$. 

$2^{\rm o}$. The integral is additive, multiplicative, and 
$\int_X F^*(x)d\phi=(\int_X F(x)d\phi)^*$.  

$3^{\rm o}$. If $F$ is a scalar function $F(x)=f(x)\cdot 1$, where $f\in C_b(X)$, then $\int_X F(x)d\phi=\phi(f)$.

$4^{\rm o}$. The integral is functorial in $X$: if there is a proper continuous map $h:Y\to X$, and $\tilde F=F\cdot h:Y\to\M(D)$, $\phi=\psi \cdot h^*: C_0(X)\to C_0(Y)\to \M(D)$, then $\int_X F(x)d\phi=\int_Y \tilde F(y) d\psi$.

$5^{\rm o}$. Suppose $U$ is an open neighborhood of the diagonal in $X\times X$, and $U_x\times\{x\}=U\cap (X\times\{x\})$. Assume that for any $x\in X$ and any $f\in C_c(U_x)$, we have: $\phi(f)F(x)\in D$. Then $\int_X F(x)d\phi=0\; \mod D_{lc}$. 

$6^{\rm o}$. $\int_X F(x)d\phi$ can be lifted to $Q_{C_0(X)}(D)$ as a properly supported element.
\end{theo}

The following simple lemma (\cite{Ka75}, section 3, lemma 1) provides the necessary estimates for the proof:

 \begin{lem} Let $B$ be a \Cst-algebra, and elements $F_1,...,F_n;\alpha_1,...,\alpha_n\in B$ satisfy the following conditions: $\sum_{i=1}^n \alpha_i^*\alpha_i=1$, and $||F_i||\le c$ for any $i$. Then $||\sum_{i=1}^n \alpha_i^*F_i\alpha_i||\le c$. 
 \end{lem}
 
 \Pf Using a faithful representation of $B$ in a Hilbert space $H$, we get for any $\xi,\eta\in H$: 
 $$|(\sum_i \alpha_i^*F_i\alpha_i(\xi),\eta)|=|\sum_i (F_i\alpha_i(\xi),\alpha_i(\eta))| \le \sum_i ||F_i|| \cdot ||\alpha_i(\xi)|| \cdot ||\alpha_i(\eta)||$$ 
 $$\le c (\sum_i||\alpha_i(\xi)||^2)^{1/2}(\sum_i||\alpha_i(\eta)||^2)^{1/2}=c||\xi|| ||\eta||. \;\cqfd $$

 \noindent\emph{Proof of the theorem.} The basic idea of the integral is the following. Let us assume that $X$ is compact. Let $B$ be the \Cst-subalgebra of $\M(D)/D$ generated by $\phi(C(X))$ and all elements $F(x),\, x\in X$. Then there is a unital homomorphism: $C(X)\to B$ which maps $C(X)$ to the center of $B$. The map $F:X\to B$ represents an element of $C(X,B)\simeq C(X)\otimes B$. Consider the multiplication homomorphism: $C(X)\otimes B\to B$. The image of $F$ under this homomorphism is by definition $\int_X F(x)d\phi\in B\subset \M(D)/D$. 
 
 Now let us translate this into the Riemann operator integration context. We will still continue to assume for the moment that $X$ is compact. The integral is constructed in the following way. Taking a finite covering $\{U_i\}$ of $X$, points $x_i\in U_i$, and a partition of unity $\sum_i \alpha_i^2(x)=1$ associated with $\{U_i\}$ (all functions $\alpha_i$ are non-negative), we consider the integral sum: $\Sigma(\{U_i\},\{\alpha_i\},\{x_i\})=\sum_i \phi(\alpha_i)F(x_i)\phi(\alpha_i)$. We assume that $||F(x)||\le c$ for all $x\in X$, so by lemma 2.5 the norm of the integral sum is $\le c$. 
 
 We will call $\{U_i\}$ an $\epsilon$-covering if for any $i$ and any $x,y\in U_i$, one has $||F(x)-F(y)||\le\epsilon$. To verify that the integral sums for two $\epsilon$-coverings $\{U_i\}$ and $\{V_j\}$ differ in norm no more than by $2\epsilon$ in $\M(D)/D$, we form the covering $W_{i,j}=U_i\cap V_j$. Let $\sum_j\beta_j^2=1$ be the partition of unity for $\{V_j\}$, and set $\gamma_{i,j}=\alpha_i\beta_j,\, x_{i,j}\in W_{i,j}$. Then we have the following estimate in $\M(D)/D$:
 $$||\sum_{i,j}\phi(\gamma_{i,j}(F(x_i)-F(x_{i,j}))\phi(\gamma_{i,j})||\le \epsilon$$
by lemma 2.5. Here 
$$\sum_{i,j}\phi(\gamma_{i,j})F(x_i)\phi(\gamma_{i,j})=\sum_i\phi(\alpha_i) F(x_i)\phi(\alpha_i)$$ 
modulo $D$ because $\sum_j \beta_j^2=1$. 
The integral is the limit in $\M(D)/D$ of the integral sums for all $\epsilon$-coverings when $\epsilon\to 0$. For a compact $X$, this gives the existence.  

In the general case of a non-compact $X$, let $\{U_i\}$ be a locally finite covering of $X$ and $\sum_i \alpha_i^2(x)=1$ the corresponding partition of unity. Then, by lemma 2.5, the norms of operators $\F_m=\sum_{i=1}^m \phi(\alpha_i) F(x_i) \phi(\alpha_i)$ are uniformly bounded in $m$. We claim that the sums $\F_m$ converge strictly in $\M(D)$ when $m\to \infty$. Indeed, for any bounded approximate unit $\{u_k\}\subset C_0(X)$ consisting of functions with compact support in $X$, our initial assumption on $\phi$ says that $\{\phi(u_k)\}$ converges strictly in $\M(D)$. Since $\F_m$ is uniformly bounded in $m$, it is enough to show that both sums $\F_m\phi(u_k)$ and $\F^*_m\phi(u_k)$ converge for any fixed $k$ when $m\to \infty$. This is true because all $u_k$ have compact support, so $\alpha_iu_k=0$ for large $i$.

Now we can take any locally finite $\epsilon$-covering of $X$ and form the integral sum $\F=\sum_{i=1}^\infty \phi(\alpha_i) F(x_i) \phi(\alpha_i)$. Modulo $D_{lc}$ any two such sums for two different $\epsilon$-coverings will differ in norm by $\le 2\epsilon$ (because all  corresponding finite sums $\F_m$ differ by $\le 2\epsilon$, as shown above). So in $Q_{C_0(X)}(D)/D_{lc}$ these sums converge in norm when $\epsilon\to 0$. The lifting of the limit to $Q_{C_0(X)}(D)$ is denoted $\int_X F(x)d\phi$. 

The proof of all properties listed in the theorem (except $5^{\rm o}$ and $6^{\rm o}$) for the case of a non-compact $X$ is the same as for a compact $X$ (see \cite{Ka75}, section 3, theorem 1). For example, for the multiplicativity property, one needs to estimate the difference modulo $D_{lc}$ of the two sums:
$$\sum_{i,j}\phi(\alpha_i)F_1(x_i)\phi(\alpha_i)\phi(\alpha_j)F_2(x_j)\phi(\alpha_j)-\sum_{i,j}\phi(\alpha_i)F_1(x_i)\phi(\alpha_i)\phi(\alpha_j)F_2(x_i)\phi(\alpha_j).$$
Both sums converge in the strict topology, so one needs only to show that all finite portions of these sums differ modulo $D_{lc}$ by less than $2\epsilon c$ if $\{U_i\}$ is an $\epsilon$-covering for $F_2$ and $||F_1(x)||\le c$ for all $x\in X$. Obviously, we can leave only those summands for which $U_i\cap U_j\ne\emptyset$. Then modulo $D_{lc}$, $||F_2(x_i)-F_2(x_j)||\le 2\epsilon$, and the estimate of the difference of those two sums $\sum_{1\le i,j\le n}$ modulo $D$ comes from lemma 2.5.

To prove property $5^{\rm o}$, we have to show that on any compact subset of $X$, the integral sums converge to $0$ in $Q_{C_0(X)}(D)/D_{lc}$. This is obvious.

Concerning property $6^{\rm o}$, we can assume that all our coverings $\{U_i\}$ consist of open sets with compact closure. This means that the operator $S=\sum_i \phi(\alpha_i) F(x_i) \phi(\alpha_i)$ is properly supported. Indeed, it is easy to check that $\supp(S)\subset \cup_i (\supp(\phi(\alpha_i))\times \supp(\phi(\alpha_i)))$ in $X\times X$. Here $\supp(\phi(\alpha_i))$ means the following: The homomorphism $\phi$ maps $C_0(X)$ onto a commutative \Cst-subalgebra in $D$. The spectrum of this commutative \Cst-algebra can be identified with a closed subset $Y\subset X$. Then $\supp(\phi(\alpha_i))$ is a compact subset in $Y$, and hence in $X$.

It is clear from the construction of the integral that once we have chosen one covering $\{U_i\}$, all subsequent coverings can be chosen as the intersections of this one with other coverings like $\{V_j\}$ above. Also the corresponding partitions of unity can be chosen as the products $\{\alpha_i\beta_j\}$. So if we fix $\{\alpha_i\}$ once and for all and change only $\{V_j,\beta_j\}$, then the resulting integral will have the form: $\sum_i \phi(\alpha_i)(\int_{U_i} F(x)d\phi)\phi(\alpha_i)$ which obviously lifts to a properly supported element. \cqfd
 
 \begin{cor} In the assumptions of the theorem, suppose that $||F(x)||\to 0$ when $x\to\infty$ in $X$. Then the sums $\F_m$ in the above proof of the theorem converge uniformly in $\M(D)$.
\end{cor}
  
     \Pf For any $\epsilon>0$, we can choose $m$ and $n$ large enough so that for any $x\in \cup_{i\in [m,n]}\supp (\alpha_i)$, we have $||F(x)||\le\epsilon$. Then $||\F_m-\F_n||\le \epsilon$ by lemma 2.5. \cqfd

  \subsection{Group averaging}
  
  We keep all assumptions of the previous subsections and assume in addition that a locally compact, second countable group $G$ acts on $X$ properly, $D$ is a $G$-algebra, and the homomorphism $\phi:C_0(X)\to \M(D)$ is $G$-equivariant. In integration over $G$, we will use the left Haar measure.
 
 \begin{lem} If $F:X\to \M(D)/D$ is $G$-equivariant, then $I(F)=\int_X F(x) d\phi$ is $G$-invariant modulo $D$ and $G$-continuous in norm (i.e. the map $G\to \M(D): g\mapsto g(I(F))$ is norm-continuous). 
 \end{lem}
 
 \Pf Because the $G$-action transforms an $\epsilon$-covering into another $\epsilon$-covering, it is clear that $I(F)$ is $G$-invariant modulo $D$.  The last assertion follows from \cite{Thom_prepr}, 1.1.4. \cqfd
 
 \medskip
 
If we want to make $I(F)$ exactly $G$-invariant we need averaging over $G$. We will adapt the averaging method of \cite{CM:L^2}, proposition 1.4, to the generality that we need.
 
\begin{prop} Let $T\in \M(D)$ be an operator with support in a set $L\times L$, where $L$ is a compact subset of $X$. Then one can define the average of $T$ over $G$, denoted $Av_G(T)$ or $\int_G g(T)dg$, as a limit of integrals $\int_C g(T)dg$, in the strict topology of $\M(D)$, over the increasing net of all compact subsets $C$ of $G$, where $dg$ denotes the Haar measure of $G$. Moreover, $||Av_G(T)||\le c ||T||$, where $c$ depends only on $L$.
 \end{prop}
 
 For the proof we need the following lemma:
 
 \begin{lem} If an operator $T\in \M(D)$ has the property that $T^*g(T)=0$ and $Tg(T^*)=0$ for any $g$ outside of a compact set $K\subset G$, then for any compact subset $C\subset G$, $(\int_C g(T)dg)^*(\int_C g(T)dg)\le |K|^2||T||^2$, where $|K|$ is the Haar measure of $K$.
 \end{lem}
 
 \Pf Consider $[(\int_C g(T)dg)^*(\int_C g(T)dg)]^n$ as the $2n$-fold integral 
 $$\int_C...\int_C g_1(T^*)g_2(T)...g_{2n-1}(T^*)g_{2n}(T)dg_1...dg_{2n}.$$
 Because of our assumption on $T$, this multiple integral actually goes over the subset of $C\times...\times C$ such that $g_i^{-1}g_{i+1}\in K$ for all $i$. So we can rewrite it as a repeated integral:
 $$\int_C\int_K...\int_K g_1(T^*h_2(Th_3(T^*h_4(T...h_{2n-1}(T^*h_{2n}(T)...)dg_1dh_2...dh_{2n},$$
where $h_{i+1}=g_i^{-1}g_{i+1}$. The latter integral is estimated as $\le |C|\cdot |K|^{2n-1}||T||^{2n}$, therefore $[(\int_C g(T)dg)^*(\int_C g(T)dg)]^n\le |C|\cdot |K|^{2n-1}||T||^{2n}$. Taking the $n$-th root of both sides and letting $n\to\infty$, we get the result. \cqfd

\medskip

\noindent\emph{Proof of the proposition.} If $a\in C_c(X)$ is a function equal to $1$ on $L$, then $T\phi(a)=T$ and $T^*\phi(a)=T^*$. Because the $G$-action on $X$ is proper, there is a compact set $K\subset G$ such that $\supp(a)$ does not intersect with any $g(L)$ for any $g\in G, \; g\notin K$. So for such $g$ we have: $T^*g(T)=T^*\phi(a)g(T)=0$, and similarly for $Tg(T^*)$. In view of the previous lemma, $||\int_C g(T)dg||\le |K| \cdot ||T||$ for any compact subset $C\subset G$, and the same for $T^*$. 

Let us denote $\int_C g(T)dg$ by $I(C)$. We need to show that for any $d\in D$, the integrals $d\cdot I(C)$ and $d\cdot I(C)^*$ converge over the net of all compact $C$. Clearly, this is true if we replace $d\in D$ with $\phi(a)$, where $a\in C_c(X)$, because if we define $C_{max}(a,T)$ as the maximal compact subset of $G$ such that $\supp(a)\cap g(L)\ne \emptyset$, then both integrals will not depend on $C$ as soon as $C_{max}(a,T)\subset C$. 

Given $d\in D$ and $\epsilon>0$, we can find $a\in C_c(X)$ such that $||d-d\cdot\phi(a)||\le \epsilon$ and $||d-\phi(a)\cdot d||\le \epsilon$. By the previous lemma, both $I(C)$ and $I(C)^*$ are bounded as functions of $C$, more precisely, $||d\cdot I(C)-\phi(a)\cdot d\cdot I(C)||\le \epsilon |K|\cdot ||T||$, and similarly for $I(C)^*$. This means that $d\cdot I(C)$ varies with $C$ no more than by $\epsilon |K|\cdot ||T||$ when $C_{max}(a,T)\subset C$. This proves the convergence. The norm estimate is also clear. \cqfd

\subsection{Operators on $C^*(G)$-modules}

We keep the assumptions of the previous subsections. In addition, we assume that the space $X$ is equipped with a Borel measure, and the group $G$ acts on $X$ properly, preserving the measure. We continue to use the left Haar measure on $G$; $\mu$ will denote the modular function of $G$. In this subsection we assume also that $X/G$ is compact. We will briefly recall a definition at the beginning of section 5 \cite{Ka16}.

\begin{defi} Let $E$ be a complex vector bundle over $X$ equipped with a Hermitian metric and an isometric action of $G$. The Hilbert module $\E$ over $C^*(G)$ is defined as the completion of $C_c(E)$ in the norm corresponding to the inner product defined below. The module structure and the inner product are defined as follows:
$$(e\cdot b)(x)=\int_G g(e)(x)\cdot b(g^{-1})\cdot \mu(g)^{-1/2}dg\in C_c(E),$$
$$(e_1,e_2)(g)=\mu(g)^{-1/2}\int_X(e_1(x),g(e_2)(x))dx\in C_c(G),$$
for $e, e_1, e_2\in C_c(E), \; b\in C_c(G)$. The scalar product under the integral is the Hermitian scalar product of $E$.
\end{defi}

Note that we will most often consider $\E$ as a Hilbert module over $C^*_r(G)$ using the natural `change of coefficients' $C^*(G)\to C^*_r(G)$.

It is easy to see that a $G$-invariant operator on $L^2(E)$ defines an operator on $\E$ (concerning boundedness, see 2.12). As explained in \cite{Ka16} (before proposition 5.5), any integral operator on $L^2(E)$ with a $G$-invariant continuous kernel and proper support defines an element of $\K(\E)$, and $\K(\E)$ is the norm-closure of such operators. In particular, by averaging over $G$ a compact operator on $L^2(E)$ with compact support on $X$, we get an element of $\K(\E)$. (An easy exercise is to show that a `rank one' operator on $\E$ is equal to $Av_G(T)$ for a rank $1$ operator $T$ on $L^2(E)$.)

Averaging is also available for operators which are not compact. One uses for this a cut-off function. We remind that a cut-off function on $X$ is a positive, continuous, compactly supported function $\gtc$ such that $Av_G(\gtc)=1$. It is easy to prove (using a cut-off function) that if a $G$-invariant, properly supported operator $F$ belongs to $\Kl(L^2(E))$ (see definition 2.4), then $F\in \K(\E)$. 

In the situation when we deal with parametrices (see section 9), and assuming that $AB-1\in \Kl(L^2(E))$, where $A$ is a $G$-invariant operator and $B$ is not, averaging $B$ over $G$ allows to assume that $B$ is also $G$-invariant (cf. \cite{CM:L^2}, proposition 1.3). This means that we can pass from parametrices in $\L(L^2(E))$ to parametrices in Hilbert modules over $C^*(G)$. 

\medskip

The following result was proved in \cite{Ka16}, proposition 5.4. 

\begin{prop} Let $A$ be an operator on $C_c(E)$ which is $L^2(E)$-bounded, $G$-invari\-ant, and properly supported. Then $A$ defines an element of $\L(\E)$ with the norm $\le {\rm const}\cdot ||A||$, where $||A||$ is the $L^2$-norm of $A$, and the constant depends only on the supports of the operators $\gtc A^*A+A^*A\gtc$ and $\gtc AA^*+AA^*\gtc$, where $\gtc$ is any cut-off function of our choice.
\end{prop}

In this subsection we will be interested only in the case of $X=G$ and $X=G/K$, where $K\subset G$ is a compact subgroup.

\begin{prop} Let $E$ be a trivial line bundle. If $X=G$, then the Hilbert module $\E$ is isomorphic to $C^*(G)$, and if $X=G/K$, where $K$ is a compact subgroup of $G$, then $\E$ is isomorphic to the Hilbert submodule of $C^*(G)$ consisting of all elements $h\in C^*(G)$ such that $ah=h$ for any $a\in K$. 

More generally, let $\gtZ$ be a finite-dimensional space with a unitary representation of $K$, and $E=G\times_K \gtZ$ over $X=G/K$. Then $\E$ is isomorphic to the Hilbert submodule of $C^*(G)\otimes \gtZ$ consisting of all elements $h\in C^*(G)\otimes \gtZ$ such that $ah=h$ for any $a\in K$ (where for $h=b\otimes z$, we have: $ah=ab\otimes a(z)$).
\end{prop}

\Pf Generally we use the left action of $G$ on all spaces, but sometimes the right action is more convenient. For the right $C^*(G)$-module structure on $C^*(G)$, obviously the formulas are:
$$(e\cdot b)(x)=\int_G e(xg)\cdot b(g^{-1})dg\in C_c(E),$$
$$(e_1,e_2)(g)=\int_Y\bar e_1(x)e_2(xg)dx\in C_c(G),$$
for $e, e_1, e_2\in C_c(E)=C_c(G), \; b\in C_c(G)$. If we replace here the right action of $G$: $e(x)\mapsto e(xg)$ with the left action: $e(x)\mapsto e(g^{-1}x)\cdot \mu(g)^{-1/2}$, we will get the same formulas as in definition 2.11. 

For the second assertion, we just replace $G$ with $K\backslash G$ and use the right $G$-action on $E$ as we did above. The last assertion is similar: put $X=K\backslash G$ and $E=(G\times\gtZ)/K$ (with the left action of $K$ on $G$), and use the right $G$-action on $E$. \cqfd

\medskip

Now we will consider $\E$ as a Hilbert module over $C^*_r(G)$.

\begin{cor} If $E=G\times_K \gtZ$, then $\K(\E)$ is a \Cst-subalgebra of $C_r^*(G)\otimes\K(\gtZ)$ consisting of all elements $s\in C_r^*(G)\otimes\K(\gtZ)$ such that $as=sa=s$ for any $a\in K$. The natural representations of $C_r^*(G)$ and $\K(\gtZ)$ in $L^2(E)$ combine together and give a representation $\K(\E)\to \L(L^2(E))$ which maps $\K(\E)$ isomorphically onto the closure of $G$-invariant integral operators in $\L(L^2(E))$ (see the discussion after definition 2.11). \cqfd
\end{cor}

\section{Pseudolocality of cosymbols}

At each point of a manifold, a cosymbol of a PDO can be considered as a multiplier of some group \Cst-algebra $C^*(G)$. In the H\" ormander calculus, $G=\R^n$ (abelian case), in more general PDO calculi $G$ may be non-abelian. \emph{Pseudolocality}, one of the main properties of a cosymbol, means that a cosymbol must commute with $C_0(G)$ modulo compact operators on $L^2(G)$.

\subsection{Abelian case}

\begin{prop} Let $f=f(\xi)$ be a bounded, differentiable function on $\R^n$ such that all its first partial derivatives vanish at infinity. Denote by $\Phi$ the operator of Fourier transform on $L^2(\R^n)$, by $f$ the operator of multiplication by the function $f$ on $L^2(\R^n)$, and by $F$ the operator $\Phi^{-1}f\Phi$. Then for any function $a=a(x)\in C_0(\R^n)$ considered as a multiplication operator on $L^2(\R^n)$, the commutator $[F,a]$ is a compact operator on $L^2(\R^n)$.

The assertion remains true in the situation with a compact parameter space $Z$. More precisely, if $L^2(\R^n)$ is replaced with $L^2(\R^n)\otimes C(Z)$ and $f(\xi)$ with $f(\xi,z)\in C_b(\R^n\times Z)$, continuous in $z$ uniformly in $\xi$, satisfying the same assumption on its first derivatives in $\xi$ (uniformly in $z\in Z$), then for any $a\in C_0(\R^n\times Z)$, the commutator $[F,a]$ belongs to $\K(L^2(\R^n))\otimes C(Z)$.
\end{prop}

\Pf We will prove the Fourier-dual assertion. Recall that $\Phi^{-1}C_0(\R^n)\Phi=C^*(\R^n)$, the \Cst-algebra of the abelian group $\R^n$. The algebra $C^*(\R^n)$ contains the dense subalgebra $C_c(\R^n)$ (compactly supported continuous functions) with convolution as multiplication. We need to prove that for any $b\in C_c(\R^n)$, the commutator $[b,f]\in \K(L^2(\R^n))$.

The assumption on the first derivatives of $f$ implies that if $\xi\to\infty$ (or $\eta\to\infty$) and $||\xi-\eta||$ remains bounded, then $|f(\xi)-f(\eta)|\le \int_0^1 |\pt f/\pt t (\xi+t(\eta-\xi))|dt\to 0$. The commutator $[b,f]$ is an integral operator with the kernel $k(\xi,\eta)=b(\xi-\eta)(f(\xi)-f(\eta))$. Since $b$ has compact support, we obviously get $\int |k(\xi,\eta)|d\xi\to 0$ when $\eta\to\infty$ and $\int |k(\xi,\eta)|d\eta\to 0$ when $\xi\to\infty$. The Schur lemma \cite{Hor_book}, 18.1.12, easily implies that the integral operator with the kernel $k$ is compact. 

The proof of the generalized version of the statement (with the parameter space $Z$) is the same. In fact, it is enough to work with $f(\xi,z)=f_1(\xi)f_2(z)$ and $a(x,z)=a_1(x)a_2(z)$. \cqfd

\begin{cor} The assertion of proposition 3.1 remains true if $f$ is bounded and measurable, but differentiable only outside of a compact subset of $\R^n$, with all first derivatives of $f$ vanishing at infinity. Moreover, the norm of the operator $F\,\mod \K(L^2(\R^n))$ does not exceed $\limsup_{\xi\to\infty}|f(\xi)|$.
\end{cor}

\Pf We can write $f=f_0+f_1$ where $f_0$ is bounded,  measurable and has compact support, and $f_1$ satisfies the assumptions of the theorem. For any $a\in C_0(\R^n)$, the products $\Phi^{-1}f_0\Phi\cdot a$ and $a\cdot\Phi^{-1}f_0\Phi$ belong to $\K(L^2(\R^n))$ by the Rellich lemma. So the assertion of the proposition remains true. For the last statement, we can take $f_0$ with as large compact support as we want. The norm of $F\,\mod \K$ depends only on $\sup |f_1|$. \cqfd

\medskip

In addition to pseudolocality of the operator $F$, we also have the following property. In the assumptions of proposition 3.1, we can consider $f$ as a distribution and denote by $\hat f$ its Fourier transform.

\begin{prop} If $f$ satisfies the assumptions of proposition 3.1, then for any $\alpha\in C_c^\infty(\R^n)$ with $\alpha(0)=0$, the operator $F_\alpha$ with the convolution kernel $\alpha\cdot  \hat f$ belongs to $C^*(\R^n)$.
\end{prop} 

\Pf Write $\alpha=\sum_{j=1}^n a_jx_j$, where $\{x_j\}$ are coordinate functions in $\R^n$ and $a_j\in C_c^\infty(\R^n)$. Let $\{\xi_j\}$ be the Fourier-dual coordinates. For any $j$ we have: $a_jx_j\hat f=(-i)a_j\cdot \widehat{\pt f/\pt \xi_j}$. Here $\pt f/\pt \xi_j\in C_0(\R^n)$, and $a_j$ can be considered as a rapidly decreasing function. Therefore $\hat a_j$ is also rapidly decreasing, and the convolution $\hat a_j * \pt f/\pt \xi_j\in C_0(\R^n)$. Passing to Fourier transforms we get $a_j\cdot \widehat{\pt f/\pt \xi_j}\in C^*(\R^n)$, which proves our assertion. \cqfd

\subsection{Non-abelian case}

The previous statements were related with the convolution operator $F$ on the $L^2$-space of the translation group $\R^n$. Now we consider more general groups. 

Let $\gtG$ be a Lie group with a fixed left-invariant Haar measure. The modular function of $\gtG$ will be denoted $\mu$. By definition, the algebra $C_r^*(\gtG)$ is the closure in $\L(L^2(\gtG))$ of the set of all compactly supported convolution operators on $\gtG$ with $L^1$-kernels. We will also consider $C_0(\gtG)$ as a subalgebra of $\L(L^2(\gtG))$.

In the following, we will need to consider linear operators $F:C_c^\infty(\gtG) \to C_c^\infty(\gtG)$ which admit an adjoint: $C_c^\infty(\gtG) \to C_c^\infty(\gtG)$ (in the $L^2(\gtG)$ sense). Any operator $F$ of this kind has a distributional kernel, which we will denote $K_F$. Symbolically we can write: $F(\phi)(g)=\int_\gtG K_F(g,h)\phi(h)dh$. The following lemma is due to L. Schwartz (cf. \cite{Yos}, theorem 6.3.2).

\begin{lem} A map $F$ (with the above properties) is $\gtG$-equivariant iff its distributional kernel is left $\gtG$-invariant: $K_F(g,h)=K_F(tg,th)$ for all $g,h,t\in \gtG$. In this case, there exists a distribution $k$ on $\gtG$ such that $K_F(g,h)=k(g^{-1}h)$, and $k$ is the convolution kernel for $F$, i.e. (symbolically) $F(\tilde\phi)(g)=\int_\gtG k(h)\phi(gh)dh$ for any $\phi\in C_c^\infty(\gtG)$ (which means that $k$ is applied to $\phi(gh)$ over the $h$ variable).
\end{lem}

We are interested in $L^2$-bounded operators on $L^2(\gtG)$ (and similar spaces) which are defined by convolution kernels with compact support. First of all, the convolution in lemma 3.4 is associated with the right $G$-action. Replacing $k(h)$ with $\tilde k(h)=k(h)\mu(h)^{-1/2}$, we get the usual left convolution action: $F(\phi)(g)=\int_\gtG \tilde k(h)\phi(h^{-1}g)dh$. Note that the modular function $\mu$ will play little role in the future discussion because we will consider only compactly supported $k$, and $\mu$ is bounded on compact subsets.

\begin{rem} {\rm When $F$ is an $L^2$-bounded operator on $L^2(\gtG)$ with a convolution kernel $k$ which has compact support, the convolution with $k$ defines a multiplier of the convolution algebra $C_c(\gtG)$. This allows to consider $F$ as a multiplier of $C_r^*(\gtG)$.}
\end{rem}

\begin{defi} A convolutional distribution $k$ on $\gtG$ will be called non-singular away from $1$ if for any function $\alpha\in C_c^\infty(\gtG)$ which equals $0$ in some open neighborhood of $1\in \gtG$, the operator $F_\alpha$ with the convolution kernel $\alpha k$ belongs to $\Kl(L^2(\gtG))$ (i.e. $F_\alpha\cdot C_0(\gtG)\subset \K(L^2(\gtG))$ and $C_0(\gtG)\cdot F_\alpha\subset \K(L^2(\gtG))$ - see definition 2.4). 
\end{defi}

In particular, if $\alpha k$ is an $L^1$-function on $\gtG$ for any $\alpha$ as in the definition, then $k$ is non-singular away from $1$. An example of a convolution $k=\hat f$ which satisfies the conditions of definition 3.6 was given in proposition 3.3 (for $\gtG=\R^n$).

\begin{prop} Let $F$ be a $\gtG$-invariant, bounded operator on $L^2(\gtG)$ with a compactly supported convolution kernel $k$. Assume that $k$ is non-singular away from $1$. Then $F$ commutes with any $f\in C_0(\gtG)$ modulo $\K(L^2(\gtG))$. 
\end{prop}

\Pf We will use theorem 2.3. Let $a,b\in C_c(\gtG)$ have disjoint supports. We need to prove that $aFb$ is a compact operator. The operator $aFb$ acts by the formula: 
$$aFb(\phi)(g)=\int_\gtG a(g)k(h)b(h^{-1}g)\phi(h^{-1}g)dh,$$
where $\phi\in L^2(\gtG)$. 

Let us define four open neighborhoods of the point $1\in \gtG$: $U_1,U_2,U_3,U_4$, such that the closure of each $U_i$ is contained in $U_{i+1}$, with the properties: for any $h\in U_2$, $a(g)b(h^{-1}g)=0$ for any $g$; $\supp(k)\subset U_3$. Let $\alpha\in C^\infty_c(\gtG)$ be a function such that $0\le \alpha\le 1$, $\alpha$ is zero inside $U_1$ and outside $U_4$, $\alpha=1$ in $U_3-U_2$.  

Then replacing $k$ with $(1-\alpha) k$ in the above displayed formula we get the operator $aF_{1-\alpha} b$ which is $0$ because $a(g)b(h^{-1}g)$ is identically $0$ inside $U_2$, $1-\alpha=0$ in $U_3-U_2$, and $k(h)=0$ for $h$ outside $U_3$. On the other hand, if we replace $k$ in the above displayed formula with $\alpha k$, the operator $aF_\alpha b$ will be compact because by our assumption $k$ is non-singular away from $1$ and $\alpha=0$ in $U_1$. \cqfd

\begin{rem} {\rm In the assumptions of proposition 3.7, cutting the support of $k$ by any compact piece which does not contain $1\in G$ will not change $F$ modulo $\Kl(L^2(\gtG))$.}
\end{rem}

Now we return to the Hilbert modules defined in subsection 2.4. We will consider the case of $X=\gtG/\gtK$, where $\gtK$ is a compact subgroup of $\gtG$, and define convolution operators on the Hilbert module $\E$ (over $C^*_r(\gtG)$) associated with the vector bundle $E=\gtG\times_\gtK \gtZ$ over $X$, where $\gtZ$ is a finite-dimensional space with a unitary representation of $\gtK$. We denote by $\tilde E$ the trivial bundle $\gtG\times \gtZ$ over $\gtG$ and by $\tilde E_c^\infty$ the space of its compactly supported smooth sections. 

Let $K_F$ be a $\gtZ$-valued distributional kernel on $\gtG$. Then the operator $F$ on $\tilde E_c^\infty$ is defined by the same formula as on $C_c^\infty(\gtG)$ before: $F(\phi)(g)=\int_\gtG K_F(g,h)\phi(h)dh$. We will assume that $K_F$ is left $\gtG$-invariant (as in lemma 3.4) and right $\gtK$-invariant (in each variable separately). Because of the $\gtK$-invariance of $K_F$, the operator $F$ descends to the space of compactly supported sections of $E$. An easy check shows that because of $\gtG$-invariance of $K_F$, the operator $F$ preserves the inner product of $\E$. Let us assume that $F$ is bounded on $L^2(E)$. Then boundedness of $F$ on $\E$ follows by proposition 2.12. 

We will interpret $\K(\E)$ as a subalgebra of $C_r^*(\gtG)\otimes \K(\gtZ)$ as in corollary 2.14. Therefore proposition 3.7 gives sufficient conditions for $F$, considered as an operator on $L^2(E)$, to commute modulo compact operators with continuous functions from $C_0(\gtG/\gtK)$.

\begin{rem} In the case when $\gtG$ is a semisimple Lie group and $\gtK$ its maximal compact subgroup, there may be a better way to deal with convolutional distributions on $\gtG/\gtK$. In the case when $E$ is just a trivial one-dimensional bundle over $\gtG/\gtK$, the algebra $\K(\E)$ coincides with the subalgebra of $\gtK$-biinvariant elements of $C^*(\gtG)$ (for $f\in C_c(\gtG): \;f(agb)=f(g)$ for any $a,b\in \gtK, g\in \gtG$) -- see 2.14. This subalgebra is commutative (see \cite{Helg}, 10.4.1), therefore it is an algebra of continuous functions on a locally compact space. This space is the space of positive definite spherical functions on $\gtG$ (see \cite{Helg}, section 10.4). There is an analog of a Fourier transform, therefore one can define symbols.
\end{rem}

\section{Manifolds with a tangent Lie structure}

 First, some notation which will be used throughout the paper. $X$ will be a complete Riemannian manifold. We will assume that $X$ is connected. We denote by $T(X)$ the tangent bundle of $X$, by $T^*(X)$ the cotangent bundle and by $p:T^*(X)\to X$ the projection. We will usually identify $T(X)$ and $T^*(X)$ via the Riemannian metric of $X$. The tangent manifold will be denoted $TX$, and the projection $TX\to X$ also by $p$. 
 
 The following abbreviation will be frequently used: the space $T_x(X)\simeq T^*(X)$ will be denoted $\tau_x$. The elements of $\tau_x$ will be usually denoted $(x,\xi)$.
 
\subsection{Main definitions}
 
 \begin{defi} We will call a tangent Lie structure on $X$ the following set of data: 

$1^{\rm o}$ A smooth bundle $\gtG(X)\to X$ of simply connected Lie groups $\{\gtG_x, x\in X\}$.

 $2^{\rm o}$ The corresponding Lie algebroid $\gtg(X)$ with fibers $\gtg_x$. (The anchor map is zero.)

$3^{\rm o}$ A smooth subbundle $\gtK(X)\subset\gtG(X)$ consisting of maximal compact subgroups $\gtK_x\subset\gtG_x$ (for $x\in X$).

$4^{\rm o}$ The corresponding Lie algebroid $\gtk(X)$ with fibers $\gtk_x$.
 
$5^{\rm o}$ A fiber-wise linear isomorphism: $\gtg(X)/\gtk(X)\simeq T(X)$ smooth in $x\in X$.
\end{defi}

Note that a `smooth bundle' of groups does not mean a principal fiber bundle. It means a Lie groupoid with multiplication given by the fiber-wise group product. The groups in the fibers $\gtG_x$ may not be all isomorphic. 

On the other hand, all groups $\gtK_x$ are isomorphic because a compact semisimple group cannot be continuously deformed into a non-isomorphic one. So we will denote a generic group $\gtK_x$ by $\gtK$. 

In the case when all algebras $\gtg_x$ are solvable, the subbundle $\gtk(X)$ is zero. In this case, we have a smooth fiber-wise linear isomorphism: $\gtg(X)\simeq T(X)$. An example of this kind arises in the well known case of filtered manifolds, where all $\gtg_x$ are nilpotent (see section 6).

We already have a Riemannian metric on $T(X)$, and we choose also a Riemannian metric on $\gtk(X)$. Considering linearly $\gtg(X)$ as an orthogonal direct sum $T(X)\oplus \gtk(X)$ gives a Riemannian metric on $\gtg(X)$ and a left-invariant Riemannian metric on $\gtG_x$ for all $x\in X$. In particular, we have: $L^2(\gtg_x/\gtk_x)\simeq L^2(\tau_x)$.

The notation $L^2(\gtG(X))$ and $L^2((\gtG/\gtK)(X))$ will be used for the Hilbert modules over $C_0(X)$ corresponding to the fields of Hilbert spaces $\{L^2(\gtG_x),x\in X\}$ and $\{L^2(\gtG_x/\gtK_x),x\in X\}$ respectively.

We will consider $C_r^*(\gtG_x)$ for any $x\in X$ as a subalgebra of $\L(L^2(\gtG_x))$. The fields of \Cst-algebras $\{C_r^*(\gtG_x), x\in X\}$ and $\{C^*(\gtG_x), x\in X\}$ are continuous fields. Basic continuous sections are convolution operators in the fibers $\gtG_x$ with continuous compactly supported kernels which vary continuously over $x\in X$. We will denote the algebra of continuous sections of these fields (vanishing at infinity of $X$) by $C_r^*(\gtG(X))$ and $C^*(\gtG(X))$ respectively. 

\begin{notat} We will use the following notation in the rest of the paper: Let $\gtZ$ be a complex vector bundle over $X$ endowed with a Hermitian metric and a fiber-wise unitary action of the group $\gtK$. We will denote by $\E_x$ the Hilbert module over $C_r^*(\gtG_x)$ defined in 2.11 for the vector bundle $E_x=\gtG_x\times_{\gtK_x} \gtZ_x$ over $\gtG_x/\gtK_x$. The natural representaion $\K(\E_x)\to\L(L^2(E_x))$ is defined in 2.14. We will denote the Hilbert module (over $C_r^*(\gtG(X))$) of continuous sections of the field of Hilbert modules $\{\E_x\}$ (vanishing at infinity of $X$) by $\E(X)$, and the corresponding algebras of the field of compact operators and their multipliers by $\K(\E(X))$ and $\L(\E(X))$ respectively.
\end{notat} 

\subsection{Exponential map}

In our coarse PDO calculus we will need an exponential map $\exp_x$. It maps a small open ball $V_x\subset\tau_x$ with center $0\in \tau_x$ diffeomorphically onto a small neighborhood $U_x$ of the point $x\in X$. This map is defined using a connection on $X$. We will assume that $X$ is equipped with a connection $\nabla$ compatible with the Riemannian metric. We will \emph{not} assume in general that the torsion of this connection is zero. 

The exponential map will be used in this paper mainly as a tool to obtain a local trivialization of the tangent bundle of the manifold. We will need the system of open balls $\{U_x,\;x\in X\}$ of the manifold $X$ which are small enough so that the corresponding exponential maps $\exp_x:V_x\to U_x\subset X$ are close to an isometry. (Here $V_x\subset \tau_x$.)

Moreover, this trivialization should be norm-continuous, i.e. when $U_x\cap U_y\ne \emptyset$, the tangent map $U_x\cap U_y\to \R^n$ corresponding to the transition $(\exp_y)_*\cdot (\exp_x)_*^{-1}:\tilde U_x\times T_x(X)\to \tilde U_y\times T_y(X)$ should be continuous in norm. (Here $n=\dim(X)$.)

Let $x\in X$ and $p\in \tau_x$. We still denote by $(\exp_x)_*$ the tangent map $T_x(X)=T_p(T_x(X))\to T_{\exp_x(p)}(X)$ for the exponential map at the point $x$ (we identify $T_x(X)$ with its tangent space at any point). Let $\gamma=\{\gamma_t=\exp(tp),\; 0\le t\le 1\}$ be the geodesic segment joining the points $x$ and $\exp_x(p)$.

To describe the tangent map $(\exp_x)_*$ more explicitly, we need to consider infinitesimal geodesic variations of the geodesic segment $\gamma$. It is known that a vector field along a geodesic segment is an infinitesimal geodesic variation iff it is a Jacobi vector  field (\cite{KN}, ch. 8, 1.2; or \cite{Mil}, 14.3). Let us denote by $\gamma'_t$ the $t$-derivative of $\gamma_t$, by $\T$ the torsion tensor, and by $\cR$ the curvature tensor for $X$. Then the differential equation for a Jacobi field $W(t)$ along $\gamma$ is
 $$\nabla_{\gamma'_t}^2(W)+\nabla_{\gamma'_t}(\T(W,\gamma'_t))+\cR(W,\gamma'_t)\gamma'_t=0$$
(see \cite{KN}, ch. 8, section 1).

Let us choose an orthonormal frame of parallel vector fields $e_1,...,e_n$ along the geodesic segment $\gamma$ and write $W(t)=\sum_{i=1}^n W_i(t)e_i(t)$. The Jacobi equation becomes:
$$ d^2W_i/dt^2+\sum_j (\T(e_j,\gamma'_t),e_i)dW_j/dt$$
$$+\sum_j ((\nabla_{\gamma'_t}\T(e_j,\gamma'_t),e_i)+(\cR(e_j,\gamma'_t)\gamma'_t,e_i))W_j=0.$$

Recall now that there is a natural `fundamental' one-form $\theta$ on $X$ (with values in $T(X)$). It is defined by the property: $\inter_\xi(\theta)=\xi$ for any section $\xi$ of $T(X)$ (see \cite{BGV}, definition 1.21). The torsion tensor $\T$ is the exterior covariant derivative of $\theta$: $\T=\nabla\theta$. Because $\nabla^2=\cR$, we have: $\nabla(\T)=\cR\wedge\theta$. Therefore, the Jacobi equation has the following form:
$$W''(t)+A(t)W'(t)+B(t)W(t)=0.$$
Here $A$ and $B$ are $n\times n$ matrices, and $W$ is column vector of dimension $n$ ($n=\dim(X)$). The matrix $A$ 
is calculated from the tensor $\T$ and $B$ - from the tensor $\cR$. 

Consider the geodesic variation of the initial geodesic segment $\gamma$ which will produce $W$. Denote by $v(s),\; -\epsilon \le s \le \epsilon$ a linear segment in $T_x(X)$ such that $v(0)=0$ and $v'(0)=q\in T_x(X)$. Put $\alpha(t,s)=\exp_x(tp+tv(s))$ (with $s$ as the variation parameter). This geodesic variation leaves the initial point $x\in \gamma$ fixed.

Then $W(0)=0$ and $W(1)=d/ds (\exp_x(p+v(s)))|_{s=0}=(\exp_x)_*(v'(s)|_{s=0})=(\exp_x)_*(q)$.  Also $W'(0)=\nabla_{\gamma'_t}(W)|_{t=0}=q$. To show this, let us denote by $\partial/\partial t$ and $\partial/\partial s$ the vector fields of the $t$ and $s$ partial derivatives for the map $(t,s)\mapsto \exp_x(tv(s))$. One has: $\nabla_{\partial/\partial t}(\partial/\partial s)-\nabla_{\partial/\partial s}(\partial/\partial t)=T(\partial/\partial t,\partial/\partial s)$ (because, obviously, $[\partial/\partial t,\partial/\partial s]=0$). Note, however, that $(\exp_x)_*$ at the point $x$ is the identity map, and $\partial/\partial s (tv(s))=0$ at $s=0, t=0$. Therefore $T(\partial/\partial t,\partial/\partial s)_{s=0, t=0}=0$. This gives: $\nabla_{\partial/\partial t}\partial/\partial s (\exp(tv(s)))|_{s=0, t=0}=\nabla_{\partial/\partial s}\partial/\partial t (\exp(tv(s)))|_{s=0, t=0}=\partial/\partial s (v(s))_{s=0}=q$.

Note that $||\gamma'_t||=||p||$. We will assume that both $||\T||$ and $||\cR||$ are bounded in some neighborhood of the geodesic $\gamma$ by some constant $c\ge 1$. Then the above formulas for the matrices $A$ and $B$ imply that $||A(t)||\le c||p||$ and $||B(t)||\le c||p||^2$ (for $0\le t\le 1$). 

\begin{prop}  Assume that $c||p||$ is small enough. Then $(\exp_x)_*$ is almost an isometry in the neighborhood $V_x\subset T_x(X)$ which satisfies this assumption on $p$. Moreover, under this assumption, the transition tangent map defined above: $(\exp_y)_*\cdot (\exp_x)_*^{-1}$, depends norm-continuously on $x$ and $y$.
\end{prop}

\Pf We can rewrite the Jacobi equation as a $2n\times 2n$ matrix equation of order $1$ by introducing a new column vector 
$C(t)=(W(t),W'(t)/||p||)$ of dimension $2n$. The equation becomes: $$C'(t)=D(t)C(t), \;\;\;\hbox{\rm where}\;\;\; D=\left(\begin{array}{cc} 0 & ||p||\\ -B/||p||& -A\end{array}\right).$$
Here $||D(t)||\le c||p||$ (for $0\le t\le 1$).

The latter equation is equivalent to the Volterra equation:
$$C(t)-C(0)=\int_0^t D(u)C(u)du.$$

Let us denote by $\D$ the operator $f\mapsto \int_0^t D(u)f(u)du$ on the space of continuous matrix functions $f$ on $[0,1]$ with values in $\C^{2n}$. It is well known from the theory of Volterra equation that for any such $f$ the series $\sum_{k=0}^\infty \D^k(f)$ converges uniformly in $t$. The sum $\sum_{k=0}^\infty \D^k(C(0))$ is the solution of our Volterra equation. 

Moreover, if $||D||\le d$ and $||f||\le a$, then the series $\sum_{k=0}^\infty \D^k(f)$ is majorated by the convergent series $ad\cdot(\sum_{k=0}^\infty d^k/k!)=ade^d$. In particular, $||C(t)-C(0)||\le ||C(0)||d(e^d-1)$ for $t\le 1$. In our case, $d=c||p||$, $C(0)=(0,q/||p||)$, so $||C(t)-C(0)||\le ||C(0)||d^2e=||q||||p||c^2e$ if $d=c||p||\le 1$. 

This implies that $||W'(t)-W'(0)||/||p||\le ||q||||p||c^2e$ for $t\le 1$, i.e. $||W'(t)-W'(0)||\le ||q||||p||^2c^2e$. Finally, using the fact that $W(1)=W(1)-W(0)=W'(t_1)$ for some $t_1\le 1$, we get $||W(1)-W'(0)||\le ||q||||p||^2c^2e$, and since $W'(0)=q$, it follows that $||W(1)-W'(0)||/||W'(0)||\le ||p||^2c^2e$.

This means that $(\exp_x)_*$, which is the map $W'(0)=q\mapsto W(1)=(\exp_x)_*(q)$, is close to an isometry if $c||p||$ is small enough. Furthermore, if there are two intersecting open balls $U_x$ and $U_y$ in $X$, both satisfying the above assumption on the parameter $c||p||$, then the transition tangent map $(\exp_y)_*\cdot (\exp_x)_*^{-1}$ defined at the beginning of this subsection, will be norm-continuous in $x$ and $y$ because the map $(\exp_x)_*$ defined by the solution of the Volterra equation depends nor-continuously on the coefficients (i.e. on $D$) and the initial value $C(0)$. The inverse of the exponential map ($(\exp_x)_*^{-1}$) is also norm-continuous because it an inverse of a map which is norm-continuous and close to an isometry. \cqfd

\medskip

{\bf Conclusion:} The tangent map $(\exp_x)_*$ is close to an isometry when $V_x$ is sufficiently small, and the smaller $V_x$ is, the closer to an isometry this map is. Also the family of first derivative linear maps $(\exp_x)_*$ is continuous in $(x,\xi)$ (uniformly in $\xi$).

\medskip

All above considerations also apply to the case of the quotient space $\gtG/\gtK$ where $\gtG$ is a simply connected Lie group and $\gtK$ is its maximal compact subgroup. (The connection must be chosen $\gtG$-invariant in this case.) The corresponding Lie algebras will be denoted (as in the previous subsection) $\gtg$ and $\gtk$ respectively. The exponential map in this case will be denoted $\Exp_x: \gtg_x/\gtk_x\to \gtG_x/\gtK_x$. This map is a diffeomorphism of a small neighborhood of $0\in \gtg_x/\gtk_x$ onto a small neighborhood of the point $(\gtK_x)\in \gtG_x/\gtK_x$. The map $\Exp_x$ depends smoothly on $x$.

\begin{rem} {\rm The actual choice of both connections used in the construction of the exponential maps $\exp$ and $\Exp$ affects the class of PDOs that we obtain in the coarse PDO calculus. For example, in the case of filtered manifolds, if we aim to obtain a coarse PDO calculus compatible with the van Erp - Yuncken one, the connections should be chosen compatible with the grading of the Lie algebras $\gtg_x$ (see \cite{vEY19}, subsection 3.2).}
\end{rem}

By our identification of $\gtg_x/\gtk_x$ with $\tau_x$, we already have a neighborhood $V_x$ in $\gtg_x/\gtk_x$. We assume in addition that this neighborhood is small enough to satisfy the above diffeomorphism assumption on $\Exp_x$. Via the map $\Exp_x$, this neighborhood $V_x$ corresponds to some neighborhood $\V_x$ in $\gtG_x/\gtK_x$. 

The $L^2$-space of $V_x$ is related with the $L^2$-spaces of $U_x$ and $\V_x$ as follows. If we choose $V_x$ as a coordinate neighborhood for $U_x$ and $\V_x$, then in each of these two cases we have the corresponding metric tensors $g_{ij}(u)=g(\pt/\pt u_i, \pt/\pt u_j)$ (in the Euclidean coordinates $\{u_i\}$ of $V_x$). The operator of multiplication by the function $\chi(u)=|\det(g_{ij})|^{-1/4}$ gives an isometric isomorphism between $L^2(V_x)$ and $L^2(U_x)$ or $L^2(\V_x)$ respectively. (Of course, the functions $\chi$ are not the same in these two cases. But in both cases $\chi(u)=1$ in the center $u=0$ of the ball $V_x$.) We will identify all three $L^2$-spaces: $L^2(V_x), L^2(U_x)$ and $L^2(\V_x)$ via these isometric isomorphisms.

\section{The coarse PDO construction}

\subsection{Cosymbols and the two main assumptions}

Let $E$ be a (complex) Hermitian vector bundle over $X$. In a classical situation, a symbol $\sigma(x,\xi)$ of order $0$ is a bounded section of the bundle $p^*(E)$ over $TX$. (Here $\xi$ is a covector at $x\in X$.) So for any $x\in X$, $\sigma(x,\cdot)$ is an element of $\L(E_x)\otimes\M(C_0(\tau_x))$. The Fourier dual element $\Phi^{-1}\sigma(x,\cdot)\Phi\in \L(E_x)\otimes\M(C^*(\tau_x))$ is the cosymbol corresponding to $\sigma$. (Here we consider $\tau_x$ as an abelian group.) A cosymbol (at $x\in X$) will be denoted $\tsig_x$. So a cosymbol in the classical case is a continuous family of elements of $\L(E_x)\otimes\M(C^*(\tau_x))$.

For a manifold $X$ with a tangent Lie structure, a cosymbol of order $\le 0$ will be a continuous family of elements $\tilde\sigma_x\in \L(\E_x)$ (see notation 4.2). Cosymbols of \emph{negative order} are those which become elements of $\K(\E(X))$ after being multiplied by any $f\in C_0(X)$. We will usually ignore negative order cosymbols, so a cosymbol may be considered (locally) as an element of the quotient $\L(\E(X))/\K(\E(X))$. 

The coarse PDO construction which is described below relies on \emph{pseudolocality} (section 3): the cosymbol $\tilde\sigma_x$ must commute with $C_0(\gtG_x)$ in the Calkin algebra $\L/\K$ of $L^2(E_x)$. We will assume that for any $x\in X$, $\tilde\sigma_x\in\L(L^2(E_x))$ satisfies the conditions of proposition 3.1 (in the abelian case) or proposition 3.7. In the course of the construction (see below) the support of each $\tsig_x$ will be cut to a compact set; this guarantees the possibility to apply proposition 3.7.

 An important second assumption is \emph{norm-continuity} of cosymbols. We will use the notation $U_x, V_x, \V_x$ introduced in section 4. Let $\tilde\sigma\in\L(\E(X))$ be a cosymbol. Let us choose a continuous family of open balls $\{U_x, x\in X\}$ which satisfy the following properties: the vector bundle $\gtZ$ is trivial over each $U_x$ (and this trivialization depends continuously on $x$), and the maps $\exp_x:V_x\to U_x$ and $\Exp_x:V_x\to \V_x$ are well defined (see section 4.2). In particular, this will mean that over $\V_x$, the restriction $L^2(\gtG_x\times_{\gtK_x}\gtZ_x)|_{\V_x}$ is isomorphic to $L^2(\V_x)\otimes \gtZ_x$ (continuously in $x\in X$).

Next, pick a function $\nu\in C_0^\infty([0,1))$ such that $0\le \nu \le 1$, $\nu(t)=1$ for $t\le 1/2$, and $\nu(t)=0$ for $t\ge 2/3$. Let us assume that $V_x\subset\tau_x$ is a Euclidean ball of radius $r_x$. Define the function $\nu_x(v)$ on $V_x$ by $\nu_x(v)=\nu(||v||/r_x)$. 

Using the isomorphisms $L^2(U_x)\simeq L^2(V_x)\simeq L^2(\V_x)$ of section 4.2, let us first transplant $\nu_x$ into $L^2(\V_x)$. Then transplant the operator $\nu_x\tsig_x\nu_x$ into $L^2(U_x)\otimes \gtZ_x\subset L^2(\gtZ)$ (recall that $\gtZ$ is trivial over $U_x$). This gives an operator on $L^2(\gtZ)$. Call this operator $F(x)$.

In the abelian case, the construction of the operator $F(x)$ is the same, except that $\V_x=V_x$ for any $x$.

\begin{defi} The cosymbol $\tilde\sigma$ will be called norm-continuous if the family $\{F(x)\in \L(L^2(\gtZ))\}$ is norm-continuous in $x\in X$ modulo $\K(L^2(\gtZ))$.
\end{defi}

If we choose another function $\nu'_x$, with a smaller support, e.g. such that $\nu'_x\nu_x=\nu'_x$, then $\nu'_x\tsig_x\nu'_x=\nu'_x(\nu_x\tsig_x\nu_x)\nu'_x$. As a product of norm-continuous functions of $x$, this will be norm-continuous.

Note that the pseudolocality condition for the operators $F(x)$ still holds in the Calkin algebra of $L^2(\gtZ)$. During the transplanting procedure, the support of the initial cosymbol has been cut by multiplication with the function called $\nu$ and also the cosymbol was (twice) conjugated by the functions called $\chi$ (see the end of section 4.2). All this does not change the cosymbol modulo compact operators.

The norm-continuity condition in the case of the H\" ormander $(1,0)$ calculus and the van Erp - Yuncken calculus will be discussed in section 6.

\subsection{The coarse PDO construction - final step}

When the family $\{F(x)\}$ satisfies the assumptions of pseudolocality (section 3) and  norm-continuous (definition 5.1), we can apply theorem 2.5. We denote by $\F$ the operator integral $\int_X F(x)d\phi$ lifted to $\L(L^2(\gtZ))$. (Here $\phi$ is the natural action of $C_0(X)$ on $L^2(\gtZ)$ by multiplication.) We choose the lift so that $\F$ has proper support -- see theorem 2.5, $6^{\rm o}$. 

Note that if we choose the functions $\nu_x$ in the construction of the family $\{F(x)\}$ with smaller supports, this will not change the operator integral - see theorem 2.5, $5^{\rm o}$. 

\medskip

Recall that negative order symbols were defined as such elements  $\tsig\in\L(\E(X))$ that $f\tsig\in \K(\E(X))$ for any $f\in C_0(X)$. \emph{We define negative order operators as elements of $\Kl(L^2(\gtZ))$} (see definition 2.4).

\begin{theo} For the cosymbols which are pseudolocal (section 3) and norm-continuous (definition 5.1), the correspondence between the cosymbol $\tilde\sigma$ and the operator $\F$ constructed out of it has the following properties (modulo negative order cosymbols and negative order operators):

$1^{\rm o}$ Composition of cosymbols $\mapsto$ composition of operators.

$2^{\rm o}$ $\tsig^* \mapsto \F^*$.

$3^{\rm o}$ If a cosymbol $\tsig_x$ is bounded in $\L(L^2(E_x))/\K(L^2(E_x))$ for all $x\in X$ by $C>0$, then the norm of the operator $\F$ in $\L(L^2(\gtZ))/\K(L^2(\gtZ))$ does not exceed $C$.
\end{theo}

\Pf All assertions follow directly from the construction and theorem 2.5.  \cqfd

\begin{rem}{\rm The fact that we ignore `negative order' cosymbols (i.e. take a quotient by $\K(\E(X))$) allows to consider homogeneous `classical' cosymbols, as well as cosymbols of van Erp - Yuncken \cite{vEY19}.}
\end{rem}

\subsection{Recovering the cosymbol of a coarse PDO}

In the usual H\" ormander PDO calculus, recovering the symbol of a properly supported PDO is an easy procedure: one applies the operator to the function $\exp(i(x,\xi))$. After that a Fourier transform gives the cosymbol.

In the coarse PDO calculus one can recover the cosymbol up to a small $\epsilon>0$ (or `up to homotopy'). For the index theory, this is enough. (Note that the construction that follows is modulo negative order operators.)

We assume that we have all the data concerning the tangent Lie structure: $\{\gtG_x,\gtK_x,\gtg_x/\gtk_x\simeq\tau_x\}$, the vector $\gtK$-bundle $\gtZ$, and the system of neighborhoods $U_x,V_x,\V_x$. 
The idea for recovering the cosymbol is to cut (for any $x\in X$) a small piece of our operator $\F\in L^2(E)$ which sits in the neighborhood $U_x$, transplant it into the neighborhood $\V_x\subset \gtG_x/\gtK_x$ and average it over $\gtG_x$. 

To implement this idea, let $\gtc_x$ be a cut-off function on $\gtG_x$ with the support in $\V_x$, i.e. a non-negative continuous function such that $\int_{\gtG_x}g(\gtc_x) dg=1$. The support of $\gtc_x$ can be as small and close to $1\in \gtG_x$ as necessary. First we transplant $\gtc_x$ into $U_x$ using the isomorphisms already used in the coarse PDO construction. The field $\{F(x)\}$ was assumed norm-continuous in $x$. Therefore, if the support of $\gtc_x$ is small, the product $\gtc_x\F$ will be close to $\gtc_x F(x)$, and therefore to $\gtc_x\tsig_x$. (We can certainly assume that $\gtc_x\nu_x=\gtc_x$.) 

Now going in the opposite direction to the steps in sections 5.1-5.2, we first consider $\gtc_x\F$ as an element of $\L(L^2(U_x)\otimes \gtZ_x)$, then transplant it to $\V_x$ and consider it as an element of $\L(L^2(\gtG_x\times_{\gtK_x}\gtZ_x)|_{\V_x})$. Averaging over $\gtG_x$ will give an element of $\L(\E_x)$. Since $\tsig_x$ was $\gtG_x$-invariant, averaging $\gtc_x\F$ over $\gtG_x$ we get an element close to $\tsig_x$.

\section{H\" ormander's and van Erp - Yuncken's calculi}

\subsection{H\" ormander's $\rho=1,\delta=0$ calculus}

In this case, the pseudolocality condition for operators of order $0$ is ensured by proposition 3.1, and also follows from the H\" ormander calculus. The norm-continuity condition comes easily from the symbol theory. Let $E$ be a (complex) Hermitian vector bundle over $X$. A symbol $\sigma(x,\xi)$ of order $0$ is a bounded section of the bundle $p^*(E)$ over $TX$. In a local trivialization of $T(X)$, the symbol satisfies the condition:
$$||\partial^\beta/\partial x^\beta \;\;\partial^\alpha/\partial \xi^\alpha \;\;\sigma(x,\xi)||\le C_{\alpha,K}(1+||\xi||)^{- |\alpha|},$$
for any compact subset $K\subset X$ and all $x\in K$, with constants $C_{\alpha,K}$ which depend on $\sigma$ and $K$. This implies norm-continuity of the symbol (uniform in $\xi$).

(Note that H\" ormander's calculus for $\delta>0$ does not have this norm-continuity property, and when $\delta=0, \rho<1$, the calculus is not diffeomorphism invariant.)

\medskip 

To change the trivialization, one can use the following

\begin{lem} Let $\sigma(x,\xi)$ be a function on $K\times \R^n$, where $K$ is a compact subset of $\R^n$, and let $\{\psi_x\}:(x,\xi)\mapsto (x,\psi_x(\xi))$ be a norm-continuous family of invertible linear maps of $K\times \R^n$ into itself.

$1^{\rm o}$ Assume that $\sigma(x,\xi)$ is continuous in $x$ uniformly in $\xi$. 

$2^{\rm o}$ Also assume that $\sigma(x,\xi)$ is differentiable in $\xi$, and for the exterior derivative $d_\xi$, there is an estimate: $||d_\xi \sigma(x,\xi)||\le C\cdot (1+||\xi||)^{-1}$ with the constant $C$ which does not depend on $(x,\xi)$. 

Then $\sigma(x,\psi_x(\xi))$ satisfies the same two conditions as $\sigma(x,\xi)$ (with a different constant $C$). Moreover, if the assumptions concerning $\sigma$ hold only outside of a compact subset in $\xi\in\R^n$, the assertion remains true outside of a (possibly larger) compact subset in $\xi\in\R^n$.
\end{lem}

\Pf The assertion about the second condition is clear. To show that the first condition is also preserved, let points $x,y\in K$ be so close to each other that $||\sigma(x,\xi)-\sigma(y,\xi)||\le \delta$ for any $\xi$ (by the assumption of continuity of $\hat\sigma$ in $x$ uniformly in $\xi$). We can write $\sigma(x,\psi_x(\xi))-\sigma(y,\psi_y(\xi))$ as a sum of two expressions: $(\sigma(x,\psi_x(\xi))-\sigma(x,\psi_y(\xi)))+(\sigma(x,\psi_y(\xi))-\sigma(y,\psi_y(\xi)))$. The norm of the second expression is estimated by $\delta$. 

The norm of the first expression is estimated as $\le ||d_\eta(\sigma(x,\eta))||\cdot ||\psi_x(\xi)-\psi_y(\xi)||$, where $\eta$ is some point on the segment joining $\psi_x(\xi)$ and $\psi_y(\xi)$, i.e. $\eta=t\psi_x(\xi)+(1-t)\psi_y(\xi)$ with $0\le t \le 1$. The first multiple is estimated as $\le C(1+||\eta||)^{-1}$. The second multiple is estimated as $\le ||\psi_x-\psi_y||\cdot ||\xi||$. Because $\psi_x$ is norm-continuous in $x$ and invertible for any $x\in K$, the map $\xi\mapsto \eta=t\psi_x(\xi)+(1-t)\psi_y(\xi)=t(\psi_x(\xi)-\psi_y(\xi))+\psi_y(\xi)$ is still an invertible linear map when $x$ and $y$ are close enough. Therefore $(1+||\eta||)^{-1}\le C_1(1+||\xi||)^{-1}$ for some constant $C_1$ depending on $x$ and $y$. This implies $||d_\eta(\sigma(x,\eta))||\cdot ||\psi_x(\xi)-\psi_y(\xi)||\le C_2(1+||\xi||)^{-1}\cdot ||\psi_x-\psi_y||\cdot ||\xi||$, thus proving the uniform continuity of $\sigma(x,\psi_x(\xi))$. \cqfd

\medskip

In the framework of the coarse PDO calculus we can now weaken the assumptions on the symbols of the H\" ormander calculus (for example, we do not need too many derivatives):

\begin{defi} Let $E$ be a (complex) vector bundle over $X$. A symbol $\sigma(x,\xi)$ of order $0$ is a bounded measurable section of the bundle $\L(p^*(E))$ over $TX$ satisfying the following  conditions:

$1^{\rm o}$ For any compact subset in $x\in X$, $\sigma(x,\xi)$ is continuous in $x$ uniformly in $\xi$ outside of a compact subset in $\xi$. 

$2^{\rm o}$ For any compact subset in $x\in X$, $\sigma(x,\xi)$ is differentiable in $\xi$ outside of a compact subset in $\xi$, and for the exterior derivative $d_\xi$, there is an estimate for any compact subset $K\subset X$: $||d_\xi \sigma(x,\xi)||\le C\cdot (1+||\xi||)^{-1}$ with the constant $C$ which depends only on $K$ and $\sigma$. 

We will say that a symbol $\sigma$ is of `negative order' if $||\sigma(x,\xi)||$ converges to $0$ uniformly in $x\in X$ on compact subsets of $X$ when $\xi\to\infty$. We will say that $\sigma$ is of `strongly negative order' if, additionally, $||\sigma(x,\xi)||$ converges to $0$ uniformly in $\xi$ when $x\to \infty$ in $X$. 
\end{defi}

\medskip

The second assumption of definition 6.2 guarantees pseudolocality of the symbol (by proposition 3.1) and independence of norm-continuity of a choice of trivialization (by lemma 6.1). The conditions of definition 6.2 are sufficient for the coarse PDO construction. The class of symbols of definition 6.2 includes the `classical' symbols (i.e. symbols homogeneous of order $0$ in the $\xi$ variable). 

Boundedness of PDOs of order $0$ in the H\" ormander calculus is proved in \cite{Hor}, 2.2.3. The fact that negative order symbols give `negative order' operators (i.e. elements of $\Kl(L^2(\gtZ))$) is proved in \cite{Hor}, 2.2.4. Strongly negative symbols give compact operators by corollary 2.7.

\subsection{van Erp - Yuncken's calculus on filtered manifolds}

For the definition of a filtered manifold, see e.g. \cite{Tan, CP1, vEY19}. A filtered manifold is a manifold $X$ equipped with a filtration of its tangent bundle $T(X)$ by smooth subbundles $\F^1\subset ... \subset \F^r=T(X)$ such that for any vector fields $v\in\F^i,\;w\in\F^j$, their commutator $[v,w]$ belongs to $\F^{i+j}$. (We set $\F^i=T(X)$ for $i\ge r$ and $\F^0=0$.) 

Let us consider the vector bundle of graded nilpotent Lie algebras $\gr(X)=\oplus_{i=1}^r \gr^i(X)$, where $\gr^i(X)=\F^i/\F^{i-1}$. For each $x\in X$, we denote the nilpotent Lie algebra in the fiber over $x$ by $\gr_x$. The Lie bracket of $\gr_x$ naturally comes from the Lie bracket of vector fields on $X$. The bundle $\gr(X)$ is a Lie algebroid over $X$ (with the zero anchor map). It is equipped with a family of dilations $\delta_\lambda$ (called the zoom action) which are endomorphisms of $\gr_x$ for any $x$. Namely, $\delta_\lambda$ acts on $\gr^i_x$ by multiplication with $\lambda^i$ .

The corresponding smooth bundle of simply connected graded nilpotent Lie groups will be denoted $\Gr(X)=\{\Gr_x, x\in X\}$. The group law in $\Gr_x$ is defined by the Baker-Campbell-Hausdorff formula. The groups $\Gr_x$ are called osculating groups. Topologically $\gr_x$ and $\Gr_x$ are diffeomorphic via the exponential map. The zoom action automorphisms on the algebras $\gr_x$ naturally induce the zoom action automorphisms on the groups $\Gr_x$. The group law in osculating groups varies smoothly in $x\in X$. However, this does not mean that all these groups $\Gr_x$ are isomorphic (nor that all Lie algebras $\gr_x$ are isomorphic) - see counterexample 2.4 in \cite{Ewe}.

In our approach, we work with complete Riemannian manifolds. In particular, all subbundles of the filtration $\{\F^i\}$ have a Euclidean metric on their fibers. This allows to take the orthogonal complement of $\F^{i-1}$ in $\F^i$ and identify this orthogonal complement with $\gr^i(X)$, which we will do. This gives a linear isomorphism: $\gr_x\simeq \tau_x$. In the terminology of \cite{vEY19}, definition 15, this is called a choice of a splitting. As soon as we get to this point, we can say that a filtered manifold has a tangent Lie structure. Therefore in the notation of section 4, we can put $\gr_x=\gtg_x$ and $\Gr_x=\gtG_x$. In the current situation we will use the notation $C^*(\Gr(X))$ instead of $C_r^*(\gtG(X))$.

The Haar measure on $\gtG\gtr_x$ coincides with the Lebesgue measure on $\gr_x$ under the above identification of $\Gr_x$ and $\gr_x$ via the exponential map, so $L^2(\Gr_x)\simeq L^2(\gr_x)\simeq L^2(\tau_x)$.

\smallskip

Let $\gtZ$ be a complex vector bundle over $X$ equipped with a Hermitian metric. The Hilbert $C^*(\Gr(X))$-module $\E(X)$ defined in notation 4.2 can be redefined in the current setting as $\gtZ\otimes_{C_0(X)} C^*(\Gr(X))$. A cosymbol of order $0$ may be considered as an element of $\L(\E(X))$. 

Now we will review some definitions and results from \cite{vEY19}. As in \cite{vEY19}, we will use a connection on $\gr(X)$ which is compatible with the grading (see remark 4.4 above), as well as with the Riemannian metric. 

Cosymbols of order $0$ are defined in \cite{vEY19}, definition 34, as the zoom-invariant elements of the algebra of \emph{compactly supported} convolutional distributions modulo smoothing convolutional distributions (for any $x\in X$). The pseudolocality condition for cosymbols of order $0$ comes from proposition 22 in \cite{vEY19} combined with our proposition 3.7. 

The norm-continuity condition for cosymbols of order $0$ follows from the theory of full symbols (sections 7 - 8 of \cite{vEY19}). Identifying $\Gr_x$ locally with $\gr_x$ via the exponential map, one can make a Fourier transform of a cosymbol $\tsig$. The resulting `full symbol' $\sigma$ will be invariant (at infinity) under the dual zoom action on $\widehat\gr_x$ (\cite{vEY19}, proposition 43). 

Corollary 45 \cite{vEY19} implies that (on compact subsets $K$ of $X$) $\pt/\pt x(\sigma(x,\eta))$ is bounded (uniformly in $\eta$), and $|\pt/\pt \eta(\sigma(x,\eta))|\le C_K(1+||\eta||)^{-1}$ for homogeneous (under the zoom action) covectors $\eta$. (In fact the right side of the last estimate reads in \cite{vEY19} as $C_K(1+||\eta||^{1/\deg(\eta)})^{-\deg(\eta)}$, which is equivalent.)

The first condition implies the norm continuity of the cosymbol $\tsig$. The second condition allows to apply lemma 6.1 to the full symbol restricted to homogeneous elements. (Note that only a change of trivialization of the tangent bundle which is compatible with the zoom action is allowed.) Therefore the norm-continuity condition for cosymbols does not depend on a trivialization. 

Summing up, for the coarse PDO construction corresponding to the van Erp - Yuncken calculus we can use the following definition:

\begin{defi} In the above notation, a cosymbol of order $0$ is an element $\tsig\in \L(\E))$ which has compact support in $\L(\E_x)$ for any $x\in X$ and satisfies the following conditions:

$1^{\rm o}$ $\tsig_x$ is norm-continuous in $x\in X$, and its Fourier transform $\sigma$ satisfies on compact subsets $K\subset X$ the estimate: $|\pt/\pt \eta(\sigma(x,\eta))|\le C_K(1+||\eta||)^{-1}$ for homogeneous (under the zoom action) covectors $\eta$.

$2^{\rm o}$ For any $x\in X$, the cosymbol $\tsig_x$ satisfies the conditions of definition 34 \cite{vEY19} for cosymbols of order $0$.

We will say that a cosymbol $\tilde\sigma$ is of `negative order' if $\tilde\sigma_x\in \L(\gtZ_x)\otimes C^*(\Gr_x)$ for all $x\in X$. We will say that $\tilde\sigma$ is of `strongly negative order' if, additionally, $||\tilde\sigma_x||$ converges to $0$ when $x\to \infty$ in $X$. 
\end{defi}

Concerning $L^2$-boundedness and compactness of operators in the van Erp - Yuncken calculus see \cite{DH22}, proposition 3.7.

\subsection{Comparison between `classical' and `coarse' calculi}

\begin{theo} Let $\tsig(x,D)$ be a bounded, properly supported pseudo-differen\-tial operator of order $0$ with the cosymbol $\tsig$ in the H\" ormander $\rho=1,\,\delta=0$ or in the van Erp - Yuncken calculus, and let $\F$ be the operator constructed out of $\tsig$ in the coarse PDO approach. Then $\tsig(x,D)-\F\in \Kl(L^2(E))$. 
\end{theo}

\Pf We need to prove that $f\cdot (\tsig(x,D)-\F)\in \K(L^2(E))$ for any $f\in C_c(X)$, so we need to prove that $\tsig(x,D)$ and $\F$ coincide modulo $\K$ on any compact piece $Y\subset X$. Let $\{U_i\}$ be a finite covering of $Y$ consisting of Riemannian balls $U_i$ and $\sum_i \alpha_i^2=1$ the corresponding partition of unity. The operator $\tsig(x,D)$ is equal to $\sum_i \alpha_i\tsig(x,D)\alpha_i$ modulo operators of lower order (pseudolocality). 

Let us compare the latter sum with the integral sum $\sum_i \alpha_i F(x_i)\alpha_i$ used in the construction of $\F$ in section 5.2. If diameters of the balls $U_i$ are small enough and $x_i$ are the centers of the balls $U_{x_i}$ (see section 5.1) then the functions $\nu_{x_i}$ used in the construction of section 5.1 are equal $1$ on $U_i$. For any $i$, both $\alpha_i F(x_i)\alpha_i$ and $\alpha_i\tsig(x,D)\alpha_i$ are PDOs of order $0$ with compactly supported distributional kernel. The difference between their cosymbols goes to $0$ uniformly in $i$ when the radii of the balls $U_i$ go to $0$ (on $Y$). Now the assertion follows from the usual norm estimate results: lemma 2.6, corollary 2.2.3 of \cite{Hor}, and proposition 3.7 of \cite{DH22}. 
\cqfd

\subsection{An example of a full PDO calculus with the groups $\gtG_x$ non-nilpotent}

Let $\gtG$ be the Lie group of upper triangular real $n\times n$-matrices with positive entries on the diagonal and $\gtg$ its Lie algebra. Denote by $\gtN$ the subgroup of triangular matrices with $1$'s on the diagonal and by $\gtD\simeq (\R_+^{\times})^n$ the diagonal subgroup of $\gtG$. The corresponding Lie subalgebras of $\gtg$ will be denoted $\gtn$ and $\gtd$ respectively. 

$\gtn$ is a graded nilpotent Lie algebra with an obvious zoom action. There is also a derivation action of $\gtd$ on $\gtn$, and $\gtg$ is isomorphic to the semidirect product of these two subalgebras. Moreover, the zoom action on $\gtn$ \emph{commutes} with the derivation action by $\gtd$. The group $\gtG$ is also a semidirect product of $\gtD$ and $\gtN$, with the action of $\gtD$ on $\gtN$ by conjugation, and this conjugation action commutes with the zoom action on $\gtN$. We have: $C^*(\gtG)\simeq C^*(\gtD,C^*(\gtN))$.

Presenting $X=\gtG$ topologically as a product $\gtN\times \gtD$, we have the corresponding PDO calculi on both multiples: the H\" ormander $(1,0)$ calculus on $\gtD$ and the van Erp - Yuncken calculus on $\gtN$. The van Erp - Yuncken calculus on $\gtN$ is preserved by the conjugation action of $\gtD$.

One can define a PDO calculus on $\gtG$ using convolution kernels given by products of convolution kernels corresponding to these two calculi on $\gtD$ and $\gtN$. This is an interesting calculus with non-trivial commutation rules.

\section{Basic $KK$-elements}

Here are some definitions from \cite{Ka88}, section 4.

Let $X$ be a complete Riemannian manifold and $\tau$ its cotangent bundle. Denote by ${\Cliff}(\tau,Q)$ the Clifford algebra bundle associated with the quadratic form $Q(v)=||v||^2$ on $\tau$. We denote by $Cl_\tau(X)$ the complexification of the algebra of continuous sections of ${\Cliff}(\tau,Q)$ over $X$, vanishing at infinity of $X$. With the sup-norm on sections, this is a \Cst-algebra. When a locally compact group $G$ acts on $X$ properly and the Riemannian metric is $G$-invariant, $\Ct(X)$ is a $G$-algebra. 

There are two canonical $K$-theory elements associated with $\Ct(X)$. The Dirac element $[d_X]\in K^0_G(\Ct(X))$ is defined as follows. Let $H=L^2(\Lambda^*(X))$ be the Hilbert space of complex-valued $L^2$-forms on $X$ graded by the even-odd form decomposition. The homomorphism $\Ct(X)\to \L(H)$ is given on (real) covector fields by the Clifford multiplication operators $v\mapsto \ext(v)+\inter(v)$. The (unbounded) operator $d_X$ is the operator of exterior derivation on $H$. The operator $D_X=d_X+d_X^*$ is essentially self-adjoint. The pair $(H, D_X(1+D_X^2)^{-1/2})$ defines the Dirac element $[d_X]$.

Another element is the local dual Dirac element. We cover $X$ with a smooth family of small balls $\{U_x\subset X\}$ (with centers $x\in X$). We define a radial covector field $\Theta_x$ on each $U_x$: at the point $y\in U_x$, $\Theta_x$ is given by $\Theta_x(y)=\rho(x,y)d_y(\rho)(x,y)/r_x$, where $\rho$ is the distance function, $d_y$ means the exterior derivative in the variable $y$, and $r_x$ is the radius of $U_x$. We assume that the radii $r_x$ vary smoothly over $X$. This family of balls actually defines an open neighborhood $U$ of the diagonal in $X\times X$, namely, $\{x\}\times U_x=U\cap(\{x\}\times X)$. We assume that both coordinate projections of the closure of $U$ into $X$ are proper maps.

We will consider $\Theta_x(y)$ as an element of the Clifford algebra fiber of $\Ct(U_x)$ at the point $y$. By definition, $\Theta_x^2-1\in C_0(U_x)\subset\Ct(U_x)$, so globally over $X$, the family of Clifford multiplications by covector fields $\Theta_x$ defines an element $[\Theta_X]\in \cR KK^G(X;C_0(X),C_0(U)\cdot C_0(X)\otimes \Ct(X))$, and consequently, an element of $\cR KK^G(X;C_0(X),C_0(X)\otimes \Ct(X))$. 

This element $[\Theta_X]$ may be considered as an element of the above group $\cR KK^G(X;C_0(X),C_0(U)\cdot C_0(X)\otimes \Ct(X))$ in two possible ways: when $C_0(X)$ acts on $C_0(X)\otimes \Ct(X)$ by multiplication over the first or the second tensor multiple. By definition, both possibilities give the same group: an element of this group is defined only by the operator  $T\in \L(\J)$, where $\J=C_0(U)\cdot C_0(X)\otimes \Ct(X)$ and $a(1-T^2), a(T-T^*)\in \K(\J)$ for any $a\in C_0(X)$. But this condition on $T$ does not depend of whether $a$ acts on the first or second tensor multiple of $C_0(X)\otimes \Ct(X)$ in view of the properness assumption on $U$.

\medskip

Let us discuss now the case of $X=\gtG/\gtK$, where $\gtG$ is a Lie group, $\gtK$ is its maximal compact subgroup. It was shown in \cite{Ka88}, 5.7, that such $\gtG/\gtK$ is a \emph{special} manifold (see \cite{Ka88}, 5.1). This means that there exists an element $[\eta_{\gtG/\gtK}]\in K_0^\gtG(\Ct(\gtG/\gtK))$ which satisfies one of the following two equivalent conditions:
$$[d_{\gtG/\gtK}]\otimes [\eta_{\gtG/\gtK}]=1_{\Ct(\gtG/\gtK)}\in KK^G(\Ct(\gtG/\gtK),\Ct(\gtG/\gtK)),$$
$$1_{C_0(\gtG/\gtK)}\otimes [\eta_{\gtG/\gtK}]=[\Theta_{\gtG/\gtK}].$$

In the case we discuss now, $[\eta_{\gtG/\gtK}]$ can be defined by a bounded covector field $\eta$ on $\gtG/\gtK$, i.e. $[\eta_{\gtG/\gtK}]= (\Ct(\gtG/\gtK), c(\eta))$, where $c$ means Clifford multiplication. We call a covector field $\eta$ \emph{special} if it defines a special element. This means in particular that $g(\eta)-\eta$ vanishes at infinity of $\gtG/\gtK$ for any $g\in \gtG$ -- the condition dictated by the definition of the equivariant $K$-theory group $K_0^\gtG(\Ct(\gtG/\gtK))$. For our applications we need only the case when $\gtG$ is \emph{simply connected}. We will sketch a construction of a special covector field $\eta$ on $\gtG/\gtK$ (cf. \cite{Ka88}, theorem 5.7).

\begin{prop} For a simply connected Lie group $\gtG$ there exists a special covector field $\eta$ on $\gtG/\gtK$. Moreover, one can always find an $\eta$ such that $\eta(y)=\Theta(x,y)$ where $x$ is the point  $(\gtK)\in \gtG$, $y$ is in a small neighborhood of it, and everywhere else on $\gtG/\gtK$, $||\eta(y)||=1$. 
\end{prop}

For any vector field $\xi$ on a Riemannian manifold $M$, with $||\xi||\le 1$, we will call the set $\{x\in M\}$ where $||\xi(x)||< 1$ the \emph{deficiency support} of $\xi$.

\medskip 

\emph{Sketch of proof.} Recall that such group $\gtG$ contains a series of normal subgroups $\{1\}=\gtN_0 \subset \gtN_1\subset ... \subset \gtN_m\subset\gtG$ so that all $\gtN_k/\gtN_{k-1}$ are Euclidean and $\gtG/\gtN_m$ is semisimple. We will use an induction on $m$ and a lemma:

\begin{lem} Assume that $\gtN$ is a solvable, simply connected normal subgroup in $\gtG$ and $\gtK$ is a compact subgroup in $\gtG$. If there is a special covector field on $\gtN$ and on $\gtG/\gtK\gtN$ (invariant under conjugation by $\gtK$ and with small deficiency supports), then there is one on $\gtG/\gtK$ (with the same properties). 
\end{lem}

\Pf The space $\gtG/\gtK$ is a principle fiber bundle over $\gtG/\gtK\gtN$ with fiber $\gtN$. Because $\eta_\gtN$ is $\gtK$-invariant, we can apply (for any $g\in \gtG$), the right translation map $\gtG/\gtK\times \gtN\to \gtG/\gtK:gk\times n \mapsto g(knk^{-1})k$ to the field $\eta_\gtN$ and obtain the field $(i_g)_*(\eta_\gtN)$ on the fiber over $g\gtK\gtN$. In fact, this gives a continuum of covector fields on the same fiber: they all differ by a shift $g\mapsto gn, n\in N$. To obtain exactly one field per fiber, we take an arbitrary smooth cross-section $s$ of the bundle $\gtG/\gtK$ over $\gtG/\gtK\gtN$. (It exists because the fiber $\gtN$ is contractible.) This allows to fix the source point $g$ of the field in the fiber $g\gtK\gtN$. The fiberwise covector field so defined will be denoted $\tilde\eta$. Because for $n\in \gtN$, $i_g$ and $i_{gn}$ differ by the left translation by $n$ on $\gtN$, and $n(\eta_\gtN)-\eta_\gtN$ vanishes at infinity of $\gtN$, the difference $g_1(\tilde\eta)-\tilde\eta$ restricted to any fiber vanishes at infinity of this fiber for any $g_1\in \gtG$.

The quotient map $\gtG/\gtK\to \gtG/\gtN\gtK$ defines a  homomorphism $\Ct(\gtG/\gtN\gtK)\to \Ct(\gtG/\gtK)$ (by lifting covectors from $\gtG/\gtN\gtK$ to $\gtG/\gtK$). The above construction gives an element $[\tilde\eta]\in(\Ct(\gtG/\gtK), c(\tilde\eta))\in KK^\gtG(\Ct(\gtG/\gtN\gtK),\Ct(\gtG/\gtK))$. The special element $[\eta_{\gtG/\gtK}]\in K_0^\gtG(\Ct(\gtG/\gtK))$ is the product of $[\eta_{\gtG/\gtN\gtK}]$ and $[\tilde\eta]$ given by the $KK$-product map: $K_0^\gtG(\Ct(\gtG/\gtN\gtK))\otimes KK^\gtG(\Ct(\gtG/\gtN\gtK),\Ct(\gtG/\gtK))\to K_0^\gtG(\Ct(\gtG/\gtK))$. 

The construction of the product involves two positive operators $M_1$ and $M_2$ (see \cite{Ka88}, 2.11) such that $M_1^2+M_2^2=1$. In our case, these two operators are just positive scalar functions on $\gtG/\gtK$, and the operator for this product is given by the Clifford multiplication with the covector field $\eta=M_1\eta_1+M_2\tilde\eta$, where $\eta_1$ is $\eta_{\gtG/\gtK\gtN}$ lifted from $\gtG/\gtK\gtN$ to $\gtG/\gtK$.

Assuming that $\eta_{\gtG/\gtK\gtN}$ has small deficiency support in $\gtG/\gtK\gtN$, near $\gtK\gtN$, we find that $||\eta_1||=1$ outside of a small neighborhood of $
\gtK\gtN\subset \gtG/\gtK$, and $\eta_1$ is always orthogonal to covetors going along fibers. On the other hand, $\tilde\eta$ always goes along fibers, and $||\tilde\eta||=1$ outside of a small subset in each fiber. This means that $||\eta||^2=M_1^2||\eta_1||^2+M_2^2||\tilde\eta||^2\le 1$. But because $\eta$ is a $KK$-product operator, it must have norm close to $1$ outside of a compact piece of $\gtG/\gtK$. Therefore we can change $\eta$ on this compact piece and also normalize it so that $||\eta||=1$ everywhere except some neighborhood of the intersection of deficiency supports of $\eta_1$ and $\tilde\eta$, in which neighborhood it will behave like the field $\Theta$ if both $\eta_{\gtG/\gtK\gtN}$ and $\eta_\gtN$ behaved so. \cqfd

\medskip

Returning to the proof of the proposition, note that a special covector field on a Euclidean space is just the radial covector field $\eta=d_x((1+||x||^2)^{1/2})$ normalized to norm 1 outside of a small neighborhood of $0$. (It is invariant under rotations.) Using the lemma, we first do the construction for the radical $\gtR=\gtN_m\subset \gtG$. In this case we can even assume at each step of induction that the quotient group $\gtN_{k+1}/\gtN_k$ is one-dimensional. When we finally come to the semisimple quotient $\gtG'=\gtG/\gtR$, we use essentially the same formula as for the Euclidean space and the fact that $\gtG'/\gtK$ has non-positive curvature (see details in \cite{Ka88}, 5.3.) \cqfd

\medskip

We assume now that $X$ is a manifold with a tangent Lie structure.

\medskip 

There exists a continuous (in $x\in X$) family of covector fields $\eta_{\gtG_x/\gtK_x}$ given by the above proposition. In fact, when we start the construction with the radical at the beginning, we can always choose a continuous family of normal subgroups of codimension $1$, and we do so on each step of the induction. For the semisimple part the construction is canonical (\cite{Ka88}, 5.3). 

\begin{notat} The fields of algebras $\{C_0(\gtG_x/\gtK_x)\}$ and $\{\Ct(\gtG_x/\gtK_x)\}$ over $X$ define $C_0(X)$-algebras which will be denoted $C_0((\gtG/\gtK)(X))$ and $\Ct((\gtG/\gtK)(X))$ respectively. The family of all special elements $\eta_{\gtG_x/\gtK_x}$ gives an element $[\eta_{(\GK)(X)}]$ of the $KK$-theory group $\cR KK(X; C_0(X),\Ct((\GK)(X)))$.
\end{notat}

\medskip

Recall now from \cite{Ka88}, 3.11, the homomorphism: 
$$j^G:KK^G(A,B)\to KK(C_r^*(G,A),C_r^*(G,B)).$$

In our case, we want to apply the family of homomorphisms $\{j^{\gtG_x}\}$ to the family of elements $\{\eta_{\gtG_x/\gtK_x}\}$ (for each $x\in X$). The result will be an element of the $\cR KK$-group which has $X$ as the parameter space, the first variable is $C_r^*(\gtG(X))$, and the second variable corresponds to the family of \Cst-algebras $\{C_r^*(\gtG_x,\Ct(\gtG_x/\gtK_x))\}$. According to \cite{Gr80} or \cite{Ka22}, the latter algebra is isomorphic to $\K(L^2(\gtG_x)\otimes Cl_{\tau_x})^{\gtK_x}$. Taking the composition with the restriction to the $\gtK_x$-invariant part of $L^2(\gtG_x)\otimes Cl_{\tau_x}$, which is isomorphic to $L^2(\gtG_x/\gtK_x)\otimes Cl_{\tau_x}$, and using the Morita equivalence between $\K(L^2(\gtG_x/\gtK_x))$ and $\C$, we obtain an element of the group $\cR KK(X; C_r^*(\gtG(X)),\Ct(X))$.

\begin{defi} Applying the family of homomorphisms $\{j^{\gtG_x}\}$ to the element $[\eta_{(\GK)(X)}]$ we get an element 
$[\eta_{C_r^*((\gtG/\gtK)(X))}]\in \cR KK(X; C_r^*(\gtG(X)),\Ct(X))$.
\end{defi}

\begin{notat} We are going to use the element $[\eta_{C_r^*((\gtG/\gtK)(X))}]$ in the index theorems of the next section. It appears that the appropriate element for this is actually the element constructed as above but out of the covector field $-\eta_{\gtG_x/\gtK_x}$ instead of $\eta_{\gtG_x/\gtK_x}$. We will denote this element $[-\eta_{C_r^*((\gtG/\gtK)(X))}]$.
\end{notat}

Our definition of Dolbeault element $[\D_X]$ is based on \cite{Ka16}, theorem 2.10. 

\begin{defi} We define the Dolbeault element $[\D_X]\in K^0(C_r^*(\gtG(X)))$ as $[-\eta_{C_r^*((\GK)(X))}]\otimes_{\Ct(X)}[d_X]$, where the $KK$-product used is:
$$KK(C_r^*(\gtG(X)),\Ct(X))\otimes K^0(\Ct(X))\to K^0(C_r^*(\gtG(X))).$$
\end{defi}

\begin{rem} {\rm Compare definition 7.4 above and definition 2.5 \cite{Ka16} of an element $[d_\xi]\in \cR KK(X;C_0(TX),\Ct(X))$. In the conventional setting of \cite{Ka16}, where $\gtG/\gtK=\R^n$, the Fourier dual of $d_\xi$ is $[-\eta_{C_r^*((\gtG/\gtK)(X))}]$. 
Theorem 2.10 of  \cite{Ka16} was stating that $[\D_X]=[d_\xi]\otimes_{\Ct(X)}[d_X]$. In the more general setting which we have here we actually replace $[d_\xi]$ with $[-\eta_{C_r^*((\GK)(X))}]$ and define $[\D_X]$ as the $KK$-product of definition 7.6.}
\end{rem}

\section{Index theorems}

This section contains index theorems for `h-elliptic' operators on complete Riemannian manifolds with a tangent Lie structure. The precise meaning of `h-ellipticity' that we use is explained in definition 8.1. We state all theorems (except 8.6) in the non-equivariant form to simplify notation. But all results of this section are true, with the same proofs, for $G$-invariant operators in the case of a proper isometric action of a second countable locally compact group $G$ on $X$ which preserves the tangent Lie structure. 

\medskip 

For the index theory we adopt the usual conventions:

$1^{\rm o}$ Self-adjoint operators of degree $1$ on $\Z_2$-graded Hilbert spaces have index in the $K^0$-groups.

$2^{\rm o}$ Self-adjoint operators on ungraded Hilbert spaces have index in the $K^1$-groups.

\subsection{Operators with index in $K^*(C_0(X))$}

\begin{defi}  A self-adjoint PDO (of order $0$) acting on sections of a vector bundle $\gtZ$ will be called h-elliptic if its cosymbol $\tsig\in \L(\E(\gtG(X)))$ satisfies the condition: $f\cdot (\tsig^2-1)\in \K(\E(\gtG(X)))$ for any $f\in C_0(X)$ (both in the graded and non-graded case).

The cosymbol $\tsig_F$ of an h-elliptic operator $F$ will be considered as an element $[\tilde\sigma_F]$ of the group $\cR KK_*(X;C_0(X),C_r^*(\gtG(X)))$. The index of the operator $F$ will be considered as an element $[F]\in K^*(C_0(X))$.
\end{defi}

\begin{defi} The Clifford cosymbol of an h-elliptic operator $F$ is defined as $[\tsig_F^{cl}]=[\tsig_F]\otimes_{C_r^*(\gtG(X))} [-\eta_{C_r^*((\GK)(X))}]\in \cR KK(X;C_0(X),\Ct(X))$.
\end{defi}

We start with the following analog of the Inverse Clifford Index Theorem 4.1 of \cite{Ka16}:

\begin{theo} In the assumptions at the beginning of the section, let $F$ be a properly supported h-elliptic operator on $X$ with the cosymbol $[\tsig_F]\in \cR KK_*(X;C_0(X),C^*(\gtG(X)))$. Then 
$$[\tsig_F^{cl}]=[\Theta_X]\otimes_{C_0(X)}[F]\in \cR KK(X;C_0(X),\Ct(X)).$$
\end{theo}

\Pf We will use a simplified version of the proof of theorem 4.1 of \cite{Ka16}. If the operator $F$ acts on sections of a vector bundle $\gtZ$, then the $KK$-product on the right hand side can be written as the pair $(\J,S)$ with $\J=C_0(U)\cdot L^2(\gtZ)\hotimes \Ct(X)$, and the operator $S$ is defined as the family $\{S_y\}$ of pseudo-differential operators (parametrized by $y\in X$): 
$$S_y=1\hotimes \Theta_x(y)+(1-\Theta_x^2(y))^{1/4}F(1-\Theta_x^2(y))^{1/4}\hotimes 1.$$
Here we consider $L^2(\gtZ)\hotimes \Ct(X)$ as the family of Hilbert spaces $L^2(\gtZ)\hotimes Cl_{\tau_y}$, where  $Cl_{\tau_y}$ is the fiber of $\Ct(X)$ over $y$.

We can always replace the neighborhood $U$ of the diagonal of $X\times X$ by a smaller neighborhood $\tilde U$ such that all $\tilde U_y=\tilde U\cap (X\times y)$ are balls in $X$ varying smoothly with $y$. For each $y\in X$, the operator $F$ restricted to $\tilde U_y$ is obtained by the coarse PDO construction (operator integration) over $\tilde U_y$. There is an obvious homotopy of this operator to the operator with the constant cosymbol $\tsig_{F,y}
$ over $\tilde U_y$. Also note that $\Theta_x(y)=-\Theta_y(x)$. Therefore the operator $S_y$ can be rewritten as
$$S_y=1\hotimes (-\Theta_y(x))+(1-\Theta_y(x)^2)^{1/4}\tsig_{F,y}(1-\Theta_y(x)^2)^{1/4}\hotimes 1.$$

Let us now compare this with the Clifford cosymbol. The Hilbert module for the product $[\tsig_F]\otimes_{C^*(\gtG(X))} [-\eta_{C^*((\GK)(X))}]$ is given by the family of Hilbert spaces $\{L^2(\gtG_y\times_{\gtK_y}\gtZ_y)\hotimes Cl_{\tau_y},\; y\in X\}$. The operator is given by:
$$S'_y=1\hotimes (-\eta_y) + (1-\eta_y^2)^{1/4}\tsig_{F,y}(1-\eta_y^2)^{1/4}\hotimes 1,$$
where $\eta_y$ is the special covector field over $\gtG_y/\gtK_y$.

We will use the homeomorphism between the neighborhood $\V_y\subset\gtG_y/\gtK_y$ and the neighborhood $\tilde U_y$ (replacing $\V_y$ by a smaller one if necessary). By the construction of the special covector field in proposition 7.1, we can always assume that in $\tilde U_y$ (or a smaller neighborhood) $\eta_y$ coincides with $\Theta_y$. The operator $\tsig_{F,y}$ also naturally restricts to a smaller neighborhood (up to operators of negative order -- see remark 3.8).

Because we can assume that $||\eta_y||=1$ outside of $\tilde U_y$, we can cut $\eta_y$ outside of $\tilde U_y$ for all $y\in X$. Then the two formulas for $S_y$ and $S'_y$ will coincide. \cqfd

\begin{theo}  In the assumptions at the beginning of the section, let $F$ be a properly supported h-elliptic operator on $X$ with the cosymbol $[\tsig_F]\in \cR KK_*(X;C_0(X),C_r^*(\gtG(X)))$. Then the formula for the index $[F]$ of this operator is 
$$[F]=[\tsig_F]\otimes_{C_r^*(\gtG(X))} [\D_X]\in K^*(C_0(X)),$$
where $[\D_X]\in K^0(C_r^*(\gtG(X)))$ is the Dolbeault element (definition 7.6).
\end{theo}

\Pf The formula follows from theorem 8.3 by applying $\otimes_{\Ct(X)}[d_X]$ to both sides of the formula in the statement of theorem 8.3 and using definitions 7.6 and 8.2, as well as the fact that $[\Theta_X]\otimes_{\Ct(X)}[d_X]=1_X$, the identity element of $\cR KK(X;C_0(X),C_0(X))$ (\cite{Ka88}, theorem 4.8).   \cqfd

\begin{rem} {\rm The reader probably noticed that the proof of the Atiyah-Singer index theorem, when this theorem is stated in the language of $K$-theory, becomes essentially a tautology.}
\end{rem}

\subsection{$G$-invariant operators with index in $K_*(C^*(G))$}

An index theory for $G$-invariant operators with index in $K_*(C^*(G))$ was introduced in \cite{Ka83}. It was used most significantly in relation with the statement of the Baum-Connes conjecture (see \cite{BCH}). It is also related with the realization of the discrete series representations. All this was described in section 5 of \cite{Ka16}. We return to this theory here because it is generalizable to operators on manifolds with a tangent Lie structure and may be useful in geometric applications. 

The basic definitions and facts concerning operators on $C^*(G)$-modules were already stated earlier in subsection 2.4. (They were used in sections 3 - 5.) In this section $G$ will be a second countable locally compact group acting properly and isometrically, with compact quotient, on a complete Riemannian manifold $X$. $\gtZ$ will be a complex vector bundle over $X$ with an isometric action of $G$. The definition of the Hilbert module $\cZ$ over $C^*(G)$ is given in 2.11 (one has to use the vector bundle $\gtZ$ instead of $E$).

In the case of a manifold with a tangent Lie structure, we assume that the $G$-action preserves this structure, and the action of $G$ on $\gtZ$ commutes with the action of the group $\gtK$ (see 4.2).

We assume that a $G$-invariant operator $F$ acts on sections of the vector bundle $\gtZ$. The index of a Fredholm operator $F$ acting on $\L(\cZ)$ will be denoted $\ind_{C^*(G)}(F)$. 

Let $\gtc\in C_c(X)$ be a cut-off function for the $G$-action on $X$ (i.e. a non-negative function such that $\int_Gg(\gtc)dg=1$). The projection $[\gtc]\in C^*(G,C_0(X))$ is defined as $[\gtc]=(\gtc\cdot g(\gtc)\cdot  \mu(g)^{-1})^{1/2}$. This projection gives an element $[\gtc]\in K_0(C^*(G,C_0(X)))$ (which does not depend of the choice of the cut-off function).

Let $[\tsig_F]\in \cR KK^G_*(X;C_0(X),C_r^*(\gtG(X)))$ be a cosymbol of an h-elliptic operator $F$ and $j^G(\tsig_F)\in KK_*(C^*(G,C_0(X)),C^*(G,C_r^*(\gtG(X))))$. Take the product $[\gtc]\otimes_{C^*(G,C_0(X))}j^G(\tsig_F)$. Denote the result $[\gtc\tsig_F]\in K^*(C^*(G,C_r^*(\gtG(X))))$.

Also note that $j^G([\D_X])\in KK(C^*(G,C_r^*(\gtG(X))),C^*(G))$.

\begin{theo} In the above assumptions on $G$ and $X$, let $F$ be a properly supported $G$-invariant h-elliptic operator on $X$ with the cosymbol $[\tsig_F]\in \cR KK^G_*(X;C_0(X),C_r^*(\gtG(X)))$. Then the formula for the index $[F]$ of this operator is 
$$\ind_{C^*(G)}(F)=[\gtc\tsig_F]\otimes_{C^*(G,C_r^*(\gtG(X)))}j^G([\D_X])\in K_*(C^*(G)).$$
\end{theo}

\Pf The proof is the same as in \cite{Ka16}, theorem 5.6. \cqfd

\section{Differential operators and parametrices}

Differential operators play a major role in any PDO calculus. But the index theory of differential operators usually depends on the index theory of operators of order $0$. In this section we will look at a standard procedure to reduce a differential operator which has a parametrix to an operator of order $0$. 

We will assume that there exists a full PDO calculus associated with the coarse PDO calculus that we use. However, we will try to carefully assess what is really necessary for index theory of differential operators. The most significant thing that we will need is an existence of a parametrix. A cosymbol of a differential operator $\D$ on a manifold $X$ is a family $\{D_x, x\in X\}$ of left-invariant differential operators on vector bundles $E_x=\gtG_x\times_{\gtK_x}\gtZ_x$ (see notation 4.2) obtained by freezing the coefficients of $\D$ at $x\in X$. The operators $D_x$ vary smoothly with $x\in X$.

\begin{defi} For any $x\in X$, a cosymbol parametrix for $D_x$ is a properly supported $\gtG_x$-invariant operator $P_x$ on $C^\infty_c(E_x)$ such that the closure of each of the three operators: $P_x,\; R_x=1-P_xD_x,\; S_x=1-D_xP_x$ belongs to $\K(\E_x)$ (cf. the explanation after 2.11), and the operators $R_x$ and $S_x$ are convolution operators with smooth kernel (`smoothing operators'). Globally over the manifold $X$, a cosymbol parametrix $P$ is a smooth (in $x$) family of point-wise cosymbol parametrices $P=\{P_x\}$, and $P,R,S\in \K(\E)$.
\end{defi}

We will follow \cite{A-Sk II, BJ, Kuc}. Recall that an unbounded operator $T$ on a Hilbert module $\E$ is called \emph{regular} if it is densely defined, its adjoint is densely defined, and its graph is \emph{orthocomplemented}. Denoting by $\Gamma=\{(h,T(h)), h\in \dom(T)\}$ the graph of $T$ and by $\Gamma^\perp=\{(T^*(h),-h), h\in \dom(T^*)\}$ its orthogonal complement, we must have: $\E\oplus\E=\Gamma\oplus\Gamma^\perp$.

For a symmetric operator $T$ on a Hilbert module $\E$, the condition to be regular is equivalent to the condition that $T^2+1$ has a bounded inverse $(T^2+1)^{-1}\in \L(\E)$ (cf. \cite{Kuc}, proposition 6). In this case the operators $(T^2+1)^{-1/2}$ and $T(T^2+1)^{-1/2}$ are also elements of $\L(\E)$. Furthermore, the operators $T\pm i$ also have bounded inverses which are elements of $\L(\E)$.

\begin{prop} Let $D$ be a cosymbol that has a cosymbol parametrix. Then the closure $\bar D$ of $D$ is a regular operator on $\E$. 
\end{prop}

 \Pf The proof of theorem 6.1 of \cite{A-Sk II} applies directly. \cqfd
 
 \medskip
 
 \begin{prop} If a cosymbol $D$ is symmetric and has a parametrix, then the `bounded cosymbol' $\tsig=\{D(D^2+1)^{-1/2}\}$ is a well defined element of $\L(\E)$, and $1-\tsig^2\in \K(\E)$.
 \end{prop}
 
 \Pf The previous proposition implies that $D$ is regular, so $\tsig$ is well defined. It remains to show that $1-\tsig^2=(D^2+1)^{-1}\in \K(\E)$. By our assumptions on the parametrix, both $DP-1$ and $P$ extend to elements of $\K(\E)$. Therefore $(D\pm i)P-1$ also belongs to $\K(\E)$. Multiplying $(D\pm i)P-1$ on the left by $(D\pm i)^{-1}$, we get: $P-(D\pm i)^{-1}\in \K(\E)$. It follows that both $(D\pm i)^{-1}$ and $(D^2+1)^{-1}$ belong to  $\K(\E)$. \cqfd
 
 \medskip
 
In order for a cosymbol to be part of an index theory, a crucial condition is \emph{pseudolocality}. If we assume the existence of a full PDO calculus, $\{\tsig_x\}$ must be pseudolocal. Another condition, norm continuity in $x$, also must be part of the full PDO calculus. 

\begin{prop} Suppose a cosymbol $D$ of a symmetric differential operator $\D$ has a parametrix at $x\in X$. In the case of the van Erp-Yuncken calculus on a filtered manifold $X$, $\{\tsig_x\}$ satisfies the pseudolocality condition. In general, on a manifold with a tangent Lie structure, if the order of $\D$ is $1$ (in the usual sense) then $\{\tsig_x\}$ satisfies the pseudolocality condition.
\end{prop}

\Pf We will follow the method of \cite{BJ} (using the formula from \cite{Ka88}, 4.2). Write $\tsig_x=2/\pi\int_0^\infty D_x
(D_x^2+\lambda^2+1)^{-1}d\lambda$. For $f\in C_c^\infty(\gtG_x)$ we get:
$$[\tsig_x,f]=2/\pi\int_0^\infty(D_x^2+\lambda^2+1)^{-1}((\lambda^2+1)[D_x,f]+D_x[D_x,f]D_x)(D_x^2+\lambda^2+1)^{-1}d\lambda$$
 We need to ensure that the integral expression is a compact operator and the integral converges in the norm topology. It is easy to see that any of the following two (equivalent) conditions is sufficient for this:
\begin{equation*}
\tag{*} |[D_x,f]|\le \const\cdot(D_x^2+1)^\alpha \;\;\hbox{\rm for some}\;\; \alpha<1/2, 
\end{equation*} 
\begin{equation*}
\tag{**} (D_x^2+1)^{-\alpha/2}|[D_x,f]|(D_x^2+1)^{-\alpha/2} \;\;\hbox{\rm is bounded for some}\;\; \alpha<1/2,
\end{equation*}
where $|T|$ means $(T^*T)^{1/2}$.

These conditions are obviously satisfied in the case when $D_x$ has order $1$ because then $[D_x,f]$ is a multiplication operator (with no derivatives).

On a filtered manifold one has to use $H$-order of differential operators (\cite{vEY19}, section 10). The commutator $[D_x,f]$ has lower $H$-order than $(D_x^2+1)^{1/2}$. We will also use the Sobolev scale (see \cite{DH22}, section 3.3). The operator $D_x^2+1$ is invertible, therefore it defines the Sobolev norms (see \cite{DH22}, 3.17). Since the $H$-order of $[D_x,f]$ is $\alpha$ times the $H$-order of $(D_x^2+1)$ for some $\alpha<1/2$, we get the inequality (**). \cqfd

\medskip

For higher order differential operators on manifolds with a tangent Lie structure, the pseudolocality for the `bounded cosymbol' $\{\tsig_x\}$ remains an open question.

\begin{defi} We will call a symmetric differential operator $\D$ h-elliptic if its cosymbol $D$ has a parametrix (definition 9.1) and the `bounded cosymbol' $\{\tsig_x\}$ (of proposition 9.3) is pseudolocal and norm continuous in $x$. 
\end{defi}

\begin{cor} The index of an h-elliptic differential operator $\D$ on a filtered manifold can be calculated by the formula of theorem 8.4, using the `bounded cosymbol' $\tsig$ of $\D$. 
\end{cor}

\end{document}